\documentclass{article}

\usepackage{arxiv}
\usepackage[utf8]{inputenc} 
\usepackage[T1]{fontenc}    
\usepackage[urlcolor=blue,anchorcolor=red, citecolor=blue, linkcolor=blue ]{hyperref} 
\usepackage{url}            
\usepackage{booktabs}       
\usepackage{amsfonts}       
\usepackage{nicefrac}       
\usepackage{microtype}      
\usepackage{lipsum}
\usepackage[english]{babel}
\usepackage{fancyhdr}
\usepackage[english]{minitoc}
\usepackage{tocloft}
\usepackage{tableof}
\usepackage{amsmath,amssymb,amsfonts,amsthm,physics,bm}
\usepackage[titletoc,title]{appendix}
\usepackage{blindtext}
\usepackage{xcolor,color}
\usepackage[utf8]{inputenc}
\usepackage{ragged2e}
\makeatletter
\newcommand{\mathleft}{\@fleqntrue\@mathmargin0pt}
\newcommand{\mathcenter}{\@fleqnfalse}
\makeatother
\usepackage{graphicx}
\usepackage{epsfig}
\usepackage{moreverb}
\usepackage{url}
\usepackage{wrapfig}
\usepackage{tikz}
\usepackage[framemethod=TikZ]{mdframed}
\usepackage{fancybox}
\usepackage[sectionbib]{chapterbib}
\usepackage{datetime}
\usepackage{cellspace}
\usepackage{subcaption}
\usepackage{listings}
\usepackage{listingsutf8}
\usepackage{footnote}
\usepackage{eurosym} 
\usepackage[np,autolanguage]{numprint}
\usepackage{indentfirst}
\usepackage{gensymb}
\usepackage{amscd}
\usepackage{dsfont}
\usepackage{tikz}
\usepackage{pstricks}
\usepackage{float}
\usepackage{chemist}
\usepackage{caption}
\usepackage[singlespacing]{setspace}
{\setlength{\baselineskip}{2.6\baselineskip}
\setlength{\parindent}{0.5cm}
\newcommand*\Laplace{\mathop{}\!\mathbin\bigtriangleup}

\title{A well-balanced FVC scheme for 2D shallow water flows on unstructured triangular meshes }

\author{Moussa Ziggaf$^{a,b,c}$\thanks{Corresponding author (moussa.ziggaf@um6p.ma) },\quad Imad Kissami $^{b}$,\quad Mohamed Boubekeur $^{c}$, \quad Fayssal Benkhaldoun $^{c}$\\[1em]
$^a$ ENSAO, LMCS, Complexe Universitaire, B.P. 669, 60000 Oujda, Morocco\\
$^b$ MSDA, Mohammed VI Polytechnic University  Lot 660, 43150 Ben Guerir, Maroc\\
$^c$ Université Sorbonne Paris Nord, LAGA, CNRS, UMR 7539, F-93430, Villetaneuse, France\\}

\begin{document}
	
\captionsetup[figure]{labelfont={bf},labelformat={default},labelsep=period,name={Fig.}}
\renewcommand{\figurename}{Fig.}

\maketitle
\noindent\rule{18.5cm}{0.6pt}
\begin{abstract}
We consider in this work the numerical resolution of a 2D shallow water system with a Coriolis effect and bottom friction stresses on unstructured meshes by a new Finite Volume Characteristics (FVC) scheme, which has been introduced in the preliminary works that will be cited below. Our main goal is to extend this approach to 2D unstructured formalism while preserving the physical and mathematical properties of the system, including the C-property.\\ 
First, we present our extension by preserving the advantages of the finite volume discretization such as conservation
property and the method of characteristics such as elimination of Riemann solvers. Afterward, an approach was
applied to the topography source term that leads to a well-balanced scheme satisfying the steady-state condition of
still water. A semi-implicit treatment will also be presented
in this study to avoid stability problems for the other source terms. Finally, the proposed finite volume method is verified on several benchmark tests and shows good
agreement with analytical solutions and experimental results; moreover, it gives a noticeable accuracy and rapidity improvement compared to the original approaches.
\end{abstract}
\noindent\rule{18.5cm}{0.6pt}

\section{Introduction}\label{sec1}
Water is a crucial issue for poverty reduction, sustainable development, and achieving the Millennium Development Goals.
However, until now, some 2.1 billion people, or 30\% of the world’s population, still do not have access to a safe water supply (\href{https://www.who.int/news/item/18-06-2019-1-in-3-people-globally-do-not-have-access-to-safe-drinking-water-unicef-who}{www.who.int/news/18-06-2019}).Climate change and uncontrollable human activities cause flood accidents to become more
frequent in recent decades. As a result, integrated water resources management is necessary and even indispensable for preventing floods and droughts. The preservation of the environment, the prevention, and control of the impact of natural risks are at the heart of
major socio-economic issues. 

In water resources management, numerical modelling is still an essential tool, and we can cite numerous numerical modelling applications of free surface flows to the management of water resources, environmental and ecosystem protection: simulation of flows due to a dam break, diversions of floods from a river to a water retention area, the change process of a river bed simulation, sediment or pollutant transport simulation in estuary and coastal environments, (see e.g. \cite{benkhaldoun2012numerical,benkhaldoun2013unstructured,chaabelasri2014well} etc. ).
In rivers, estuaries and coastal areas, flows are characterized by: great topographical and morphological complexity, strong affection, or even pure advection in the case of a dam break on a flat and slippery bottom (without friction), a variable space scale (starting from a dozen to a few thousand meters) and in time scale (starting from a few minutes to several days).

Consequently, when developing a numerical approach to solve a free surface flow or other models, one encounters major difficulties which result from the physical complexity of the area and numerical calculations. We consider in this paper a derivative as a formal first-order approximation of the three-dimensional free surface incompressible Navier–Stokes equations, using the so-called shallow-water assumption \cite{ferrari2004new}.
The difficulty in defining  accurate numerical schemes for such hyperbolic systems is related to their non-linear behaviour or, more generally, mathematical structure and the physical phenomena they generate. In particular, the presence of a shock front essentially causes numerical oscillations or artificial scattering, which are due to the treatment of the advection terms in the equations governing the transport of the water mass by a standard method of approximation. Another fundamental point is to get schemes that satisfy the preservation of steady states, such as still water equilibrium in the context of the shallow water system. Different approaches to satisfy the well-balanced property have been proposed (see, e.g. \cite{benkhaldoun2010new,bermudez1998upwind, bermudez1994upwind,leveque1998balancing,rebollo2003entropy}), and recent extensions to other types of homogeneous solvers can be found in \cite{castro2017well,michel2017well}. The extension of ENO and WENO schemes to shallow water equations has been studied in \cite{vukovic2002eno}. Unfortunately, most ENO and WENO schemes that solve real flows correctly are still very computationally expensive. On the other hand, numerical methods based on kinetic reconstructions have been studied in \cite{perthame2001kinetic}, but the complexity of these methods is relevant. However, most of the above mentioned works, even if they are unstructured two-dimensional methods they lead to rather complex and time consuming algorithms. Other approaches are more efficient, but to our knowledge, they are limited to one-dimensional problems or to two-dimensional Cartesian meshes.

Our main objective in the present study is to develop a class of Eulerian-Lagrangian methods and to accurately solve shallow water equations in 2D-unstructured meshes without relying on Riemann solvers. The proposed FVC scheme belongs to the class of methods that employ only physical fluxes and averaged states in their formulations. It can be interpreted as a predictor–corrector scheme. In the corrector stage, the considered equations are integrated over an Eulerian time–space control volume, whereas, in the predictor stage, the conservation laws are rewritten in an advective form and integrated along the characteristics defined by the  advection velocity field. This approach has shown its effectiveness through several test cases represented in the works, \cite{benkhaldoun2010simple,benkhaldoun2015projection,benkhaldoun2015family}, but the authors of these three papers remained limited in the Cartesian mesh formalism but, as I have mentioned before, the real problems are characterised by a great topographical and geometrical complexity, hence, the limit of this formalism. A more thorough study of the accuracy of this finite volume discretization method on unstructured meshes was always the objective in these works (see the conclusions of \cite{benkhaldoun2015projection,benkhaldoun2015family}). Therefore, the unstructured finite volume method doesn't only ensure the conservation of mass, which is an important property in the computation of fluid flows but also allows the complex geometry of the computational domain to be perfectly taken into account. For these reasons, we propose an extension of this scheme in an unstructured mesh.

Another strong point of this discretization method is that the Jacobian matrix of the system that caused the slowness of many approximation schemes is not involved in the calculation. Note also that many approximation schemes in the framework of conservation laws require a solver for the Riemann problem at each time step to reconstruct the numerical flow, which is completely avoided in our FVC scheme. Simply this approach is based on a combination of the characteristics method and a finite volume method. Unlike traditional finite volume methods, this technique integrates our equations along the characteristics curves so that numerical flux are easily calculated. 
We noticed that this approach has several advantages over others conservation law solution techniques. Indeed,
the main features of such a finite volume Eulerian-Lagrangian scheme are on the one hand, the capability to satisfy the conservation property resulting in numerical solutions free from spurious oscillations, and on the other hand, the achievement of strong stability and high accuracy for numerical solutions containing shocks or discontinuities.

We first present a projection of the shallow water system in order to find a form of a transport equation on which we will apply the techniques of the characteristics method for the purpose of evaluating our unknowns at the interfaces of our unstructured mesh. We also present original developments to take into account boundary conditions, to reduce the diffusion of the scheme or to increase the efficiency of the method.

Second, we announce a new reformulation of the FVC scheme adapted to non-uniform triangular mesh "unstructured mesh" as well as the discretization of the flow gradients and source terms while keeping the equilibrium and the C-property. It will also be seen that the proposed scheme has the ability to handle calculations of slowly varying flows as well as rapidly varying flows over continuous and discontinuous bottom beds.

This article is organized as follows: Section \ref{sec2} will be devoted to the presentation of the mathematical model and its projected velocity system, as well as its mathematical and physical properties. In Section \ref{sec3}, numerical methods are formulated for the reconstruction of the finite volume characteristics scheme on an unstructured meshes and the approximation of the source terms keeping  the scheme Well-balanced. Section \ref{sec4} is devoted to the numerical results of several test examples and some results interpretations. It is shown that our new approach achieves the expected accuracy and robustness. Section \ref{sec5} contains concluding remarks and an outlook.
\section{Mathematical Model}\label{sec2}
\subsection{Governing equation}
\begin{figure}[H]
	\centering
	\includegraphics[width=0.55\linewidth, height=0.25\textheight]{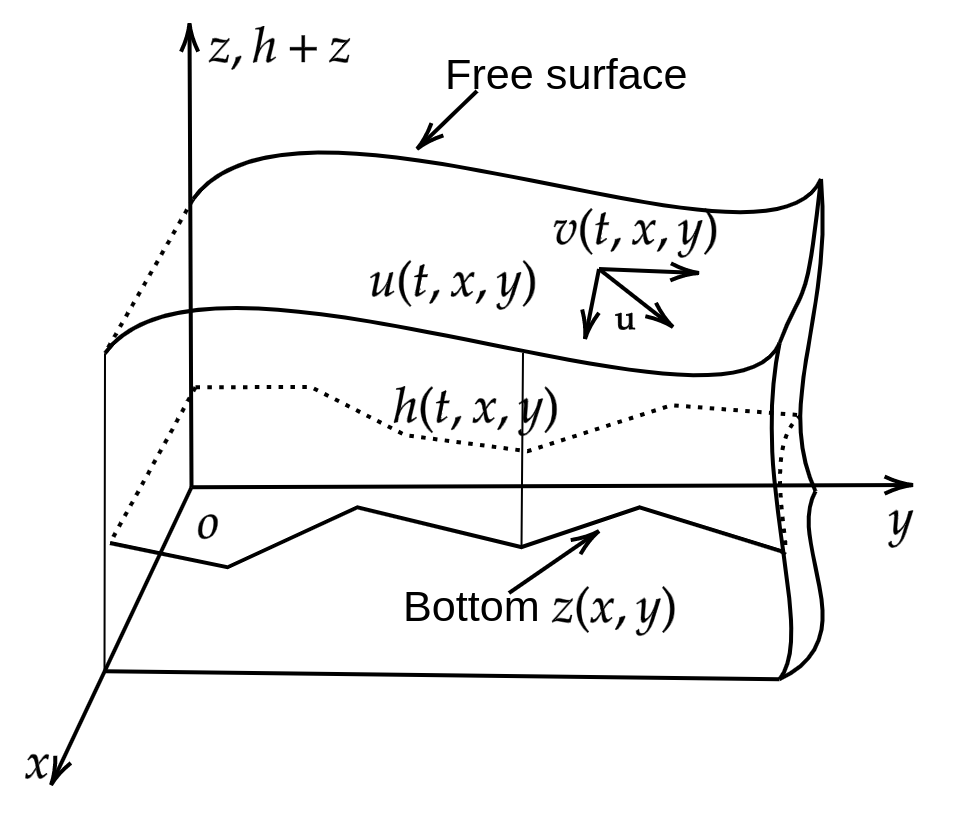}
	\caption{Illustration of shallow water model variables.}
	\label{fig:swdomain}
\end{figure}

The 2D shallow water system for the free-surface flow  with the Coriolis effect and bottom friction stress is formulated as

\begin{equation}\label{1}
	\left\{ \begin{array}{ll}
	\displaystyle  \partial_{t} h + \nabla \cdot (h\textbf{u}) = 0 \\[.3em]
	\displaystyle  \partial_{t} h\textbf{u} + \nabla \cdot (h\textbf{u} \otimes\textbf{u}) + \frac{1}{2}\grad{(gh^2)} = -gh\grad{Z} - f_c \times h\textbf{u} - r(h,\textbf{u}),  \\
\end{array}
\right.
\end{equation}
where the unknowns are always the water height $h(t,x,y) \geqslant 0 $ and  the horizontal speed
mean $\textbf{u}(t, x, y)=(u,v)^T(t,x,y) \ \in \mathbb{R}^2$. 
 The parameter $ f_c $ is linked to the angular speed of the earth's rotation, $g$ is the gravitational acceleration, $r(h,\textbf{u})$ has various expressions, for example, the asymptotic derivation mentioned in \cite{gerbeau2000derivation} leads the authors to consider, at first order, a linear friction term. The quadratic form in the Manning-Strickler velocity is nevertheless the most widely used in river flow applications \cite{delestre2010simulation,hervouet2007hydrodynamics}, 
so in this study we use the latter approximation such that the bottom's friction term $  r(h,\textbf{u}) $ is given by, $  r(h,\textbf{u}) =(r_{f_x}, r_{f_y} ) := \eta^2gh^{-1/3}|\textbf{u}|\textbf{u}$, such as $ \eta$ is the Manning roughness. The function $ Z(x, y)$ represents the bottom profile, see the \figurename{\ref{fig:swdomain}}.
In order to give the reader a global view on the  shallow water system we propose to add other aspects in the right hand side of the second equation of system (\ref{1}). For example we can add the wind's effect on the free surface, i.e. in the case where the wind is moving with high speed there is the friction term $\tau (h,\breve{\textbf{u}})$ which is not negligible. The viscosity or diffusion term can also be added if we want to solve turbulence problems in free surface flow (see e.g \cite{zaouali2008structure}). In this case, the second equation of system (\ref{1}) becomes

\begin{equation}\label{2}
\displaystyle  \partial_{t} h\textbf{u} + \nabla \cdot (h\textbf{u} \otimes\textbf{u}) + \frac{1}{2}\grad{(gh^2)} = -gh\grad{Z} - f_c \times h\textbf{u} - r(h,\textbf{u}) + \tau (h,\breve{\textbf{u}}) +\nu \overrightarrow{ \Laplace } h\textbf{u},
\end{equation}

where, $\tau(h,\breve{\textbf{u}}) = (\tau_{s_x},\tau_{s_y}) :=  \frac{1}{2}\breve{C_f}|\breve{\textbf{u}}|\breve{\textbf{u}}$, in which $\breve{\textbf{u}} = (\breve{u},\breve{v})^T$ represents the wind speed and $\breve{C_f}$  is the coefficient of wind friction with water. $\nu$ is the diffusion coefficient associated to the term,  $ \overrightarrow{ \Laplace } h\textbf{u} := (\Laplace (hu), \Laplace (hv))^T $. 

For simplicity, in this study we will not deal with these two terms ($\tau (h,\breve{\textbf{u}}) $ and $\nu \overrightarrow{ \Laplace } h\textbf{u}$ ), we rewrite the system (\ref{1})  in a conservative vectorial form
\begin{equation}\label{3}
\displaystyle  \partial_{t} W + \nabla \cdot \mathbb{F}(W) = S(W) + Q(W) ,
\end{equation}\\
where, \quad $\mathbb{F}(W) = (F(W),\ G(W))^T$.  Such that, we note $A^T$ is the transpose of a matrix $A$,   \\[2em]
\small
$W= \left( \begin{array}{c} h\\[0.5em] hu\\[0.5em] hv\\ \end{array}\right)$, \ \
$F(W)= \left( \begin{array}{c} hu\\[0.5em] hu^2+\frac{1}{2}gh^2\\[0.5em] huv\\ \end{array}\right)$, \ \
$G(W)= \left( \begin{array}{c} hv\\[0.5em] huv\\ [0.5em] hv^2+\frac{1}{2}gh^2\\ \end{array}\right)$, \ \
$S(W)= \left( \begin{array}{c} 0\\[0.5em]-gh\partial_{x} Z\\ [0.5em] -gh\partial_{y} Z\\ \end{array}\right)$, \ \
$Q(W)= \left( \begin{array}{c} 0\\[0.5em] f_chv - r_{f_x}\\[0.5em] -f_chu - r_{f_y}\\ \end{array}\right).$\\
\normalsize

Note that the equation (\ref{3}) has to be solved in a bounded spatial domain $\Omega$, with given boundary and initial conditions. In practice, these conditions depend on the phenomenon studied ( see the Section \ref{sec4} where numerical examples are discussed).

\subsubsection{Properties of the system}\label{subsec211}

A simple calculation shows that the system is still written in a quasi-linear form

\begin{equation}\label{4}
\displaystyle  \partial_{t} W + \textbf{J}_F\partial_{x}W + \textbf{J}_G\partial_{y}W  = S(W) + Q(W) ,
\end{equation}\\
where $\textbf{J}_F$ and $\textbf{J}_G$ are the Jacobian matrices of the fluxes\\
\begin{equation*}
\textbf{J}_F= \left( \begin{array}{ccc} 0&1&0\\[0.5em] gh-u^2&2u&0\\[0.5em] -uv&v&u\\ \end{array}\right)\qand
\textbf{J}_G= \left( \begin{array}{ccc} 0&0&1\\[0.5em] -uv&2v&u\\[0.5em] gh-v^2&0&2v\\ \end{array}\right).
\end{equation*}

Following the usual techniques  \cite{godlewski2013numerical, serre1999systems}, we define $\textbf{J}_{\mathbb{F}}((\alpha_1,\alpha_2)) := \alpha_1\textbf{J}_F +  \alpha_2\textbf{J}_G$. For any $(\alpha_1,\alpha_2)\in \mathbb{R}^2$, the matrix $\textbf{J}_{\mathbb{F}}((\alpha_1,\alpha_2))$ has three eigenvalues defined by
\mathleft
\begin{align*}
\lambda_1& = \alpha_1 u +\alpha_2 v, \\
\lambda_2& = \alpha_1 u +\alpha_2 v + |(\alpha_1,\alpha_2)|\sqrt{gh} ,\\
\lambda_3& =  \alpha_1 u +\alpha_2 v - |(\alpha_1,\alpha_2)|\sqrt{gh}.
\end{align*}
\mathcenter
The shallow water system is a first-order hyperbolic system of balance laws, and it is also strictly hyperbolic for $h > 0$ with real and distinct eigenvalues. \\

\textbullet) \textit{ Equilibrium }\\

An important property is related to the source terms, and the most studied balance family is related to the presence of topography's source term: the shallow water system admits non-trivial steady-states. They are characterized by\\

\begin{equation}\label{5}
\displaystyle \nabla \cdot \textbf{hu}=0, \ \ \ \  \grad(\frac{|\textbf{u}^2|}{2} + g(h+Z)) - \textbf{u}\grad\times\textbf{u}=0.
\end{equation}
For flows in complex geometry, it seems very difficult to numerically preserve all two-dimensional balances, except those that correspond to an area at rest and whose
characterization is independent of the dimension considered.\\
\begin{equation}\label{6}
\displaystyle h + Z = Cst, \ \  \ \ \ \  \textbf{u} = 0. 
\end{equation}
This particular stationary state, known as the resting lake state, is important because many flows in lakes or coastal bays are perturbations  around this balance. Therefore it is  essential to prevent numerical anomalies from disrupting the approached solution. However, the preservation of the stationary states at the numerical level is not obvious to be achieved, and even the simplest one is not an exception. In fact, (\ref{6}) correspond to a balance between flow terms and source terms, whose discretization are not correlated.

There are other categories of stationary states resulting from an equilibrium between the pressure term and the Coriolis term, i.e.

\begin{equation*}
g\grad h + f_c \times \textbf{u} = 0 .
\end{equation*}
 
For example, at large scales, the atmospheric and oceanic flows bear most of the time the perturbations of this stationary state \cite{olbers2012ocean}, therefore it is also very important to be represented in the approached solution. This balance presents a complexity added to the balance of the lake at rest because it involves non-zero speeds.\\

\textbullet) \textit{ Entropy inequality }\\

The physical system we are dealing with here is a system of conservation laws, so the energy aspect is very important in this type of system, thus we propose to say a few words on this point in order to provide insight into the treatment of this notion when constructing approximate solutions under the constraint of preserving certain properties, for example, the decay of energy in the presence of friction source terms.

We call entropy  solution to the shallow water system (see \cite{gassner2016well}), a weak solution which satisfies the following  entropy inequality\\
\mathleft
\begin{equation}\label{energie}
\partial_{t} E + \div (E + \frac{|\textbf{u}|^2}{2}) \leqslant 0,
\end{equation} 

where $E$ is a mathematical entropy (which is the mechanical energy see \cite{lions1996existence}), defined as\ \ \ 
$\displaystyle E(h,\textbf{u},Z) = h\frac{|\textbf{u}|^2}{2} + \frac{gh^2}{2} +ghZ.$

This inequality becomes equality for regular solutions, in the absence of energy loss terms, notably friction, and remains inequality, for admissible discontinuous solutions, resulting from classical calculations; we invite you to read subsection 1.1 in \cite{fjordholm2011well}. The mechanical energy, which is easily verified as convex with respect to the conservative variables, thus acts here as a mathematical entropy. In the case of the system without source terms and in the 1D problem, the mechanical energy is only one of the entropies that must be associated with the system for the problem to be properly posed.
In 2D, or when the system contains the source terms, there is no longer a complete family of mathematical entropies. Therefore, even if the inequality (\ref{energie}) alone is not sufficient for a rigorous mathematical study, it nevertheless ensures the presence of an additional bound on a certain positive function of the system unknowns and can provide information on the choice of a physical solution. The finite volume scheme presented in section \ref{sec3} verifies innately the conservation properties. The energy decay property remains more difficult to satisfy. 

\mathcenter
\subsubsection{The initial  and the boundary conditions}
To obtain a well-posed problem, we add to this system some initial condition and boundary conditions.
In this paper, we consider only two types of boundaries: solid walls on which we prescribe a slip condition and fluid limits on which we prescribe one or two conditions depending on the type of flow (subcritical, supercritical). In some cases we use the injection and Neumann homogeneous boundary conditions.\\

\begin{tabular}{cccccc}
	\hline 
	slip wall & no-slip wall &  \hspace{-1cm} subcritical inflow  &  \hspace{-1cm} subcritical outflow  & \hspace{-1cm}  supercritical inflow & \hspace{-0.1cm} supercritical outflow \\[0.5em] 
	$h_r = h_l$ & $h_r = h_l$ & $h_r =h_b$  &$h_r =h_b$  & $h_r = h_b$  & $h_r = h_l$\\ 
	 $u_r = 0 $ &$u_r = 0 $ &$u_r = u_l + \sqrt{g}(\sqrt{h_l}-\sqrt{h_r})$  &$u_r = u_l + \sqrt{g}(\sqrt{h_l}-\sqrt{h_r})$  & $u_r = u_b$  &$u_r = u_l$\\
	$v_r=v_l$  & $v_r=0$&  $v_r=0$  &$v_r=v_l$ &$v_r = 0$ &$v_r = v_l$ \\
	\hline 
\end{tabular}\\

Here the subscripts $r$ and $l$ denote the right and left  states respectively at a boundary cell interface, such as the local values of the Froude number is used to determinate whether the flow is subcritical or supercritical at a given time.
The subscript $b$ above denotes prescribed physical boundary values.
\section{Finite volume discretization of the model in 2D unstructured formalism }\label{sec3}
In this section we present the FVC scheme for the  discretization of the shallow water system (\ref{3}). The method consists of two steps and can be interpreted as a predictor-corrector approach. The first step deals with the classical finite volume method, whereas in the second step, the reconstruction of the numerical flux by constructing the intermediate state using the characteristics method \cite{ziggaf2020fvc}.
\subsection{Discretization}\label{subsec31}
\begin{figure}[H]
	\centering
	\includegraphics[width=0.95\linewidth, height=0.26\textheight]{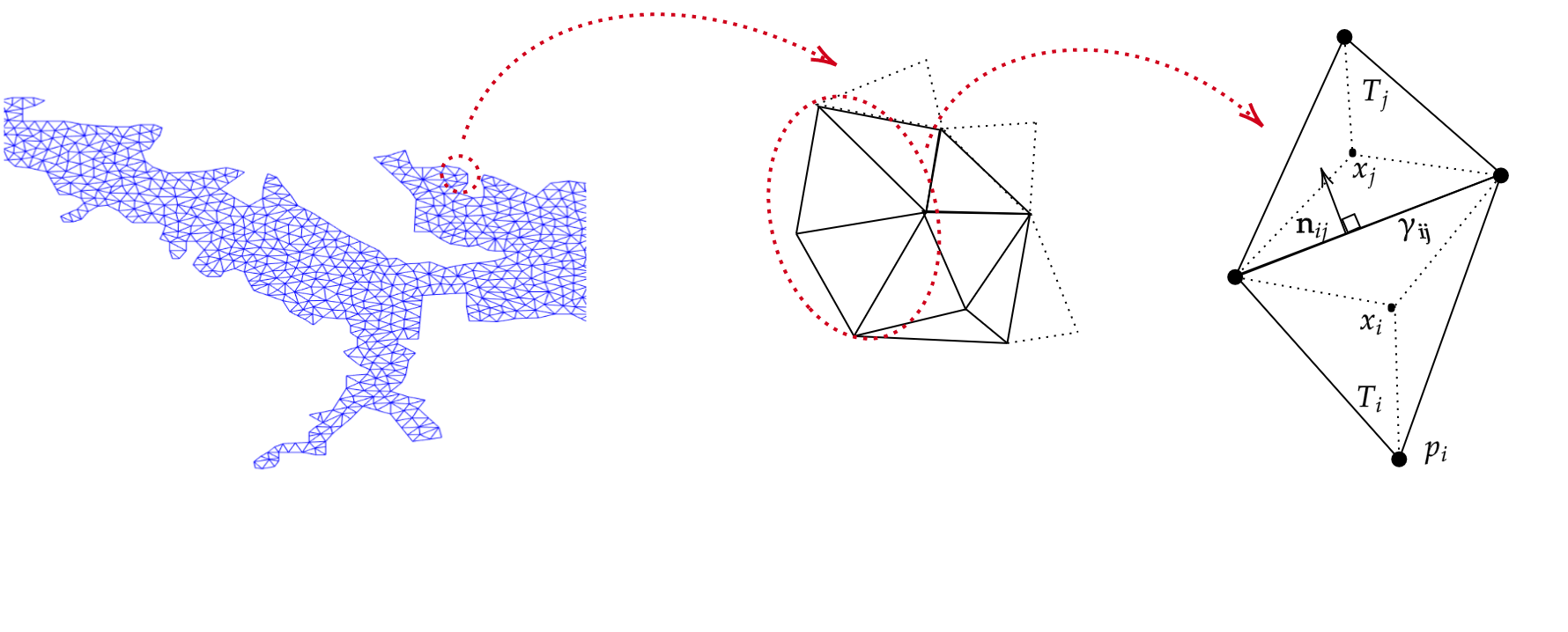}
	\caption{Generic definition of the $\Omega$ domain and the control cells of the mesh.}
	\label{fig:mesh}
\end{figure}

The integral form of the system (\ref{3}) can be written as
\begin{equation}\label{8}
\displaystyle  \frac{\partial}{\partial t} \int_{\Omega} W dV + \oint_{\partial \Omega} \mathbb{F}(W)\cdot\textbf{n} d \sigma = \int_{\Omega} ( S(W) + Q(W) )  dV,
\end{equation}

where $\Omega$ is the domain of interest, $\partial \Omega$ is the boundary surrounding, $\textbf{n}$ is the normal vector to $\partial \Omega$ in the outward direction, $dV$ and $d\sigma$ are respectively the surface element and the length element. The problem domain is first discretized into a set of triangular cells forming an unstructured
computational mesh see \figurename{\ref{fig:mesh}}. The average of conserved variables is stored at the centre of each cell, and the edges of each cell define the faces of a cell control what is called "control volume". \\

We use the following notations:\\
\begin{itemize}
	\item $x_i$, centroid of the cell $T_i$,
	\item $p_i$, vertex of  $T_i$,
	\item $\gamma_{ij}$, boundary edge between the cells $T_i$ and $T_j$,
	\item $|\gamma_{ij}|$, length of  $\gamma_{ij}$,
	\item $|T_i|$, area of the cell $T_i$,
	\item $\textbf{n}_{ij}$, unit normal to $\gamma_{ij}$, outward to  $T_i$ such as, $\textbf{n}_{ji} = -\textbf{n}_{ij}$.
\end{itemize}
For each triangular control volume, the system (\ref{3}) is written as
\begin{equation}\label{9}
|T_i|\frac{dW_i}{dt}+ \oint_{\gamma_{i}} \mathbb{F}(W)\cdot\textbf{n} d \sigma = \int_{T_i} (S(W) + Q(W) ) dV,
\end{equation}\\
where $W_i$ is the average quantity of cell $T_i$ stored at the cell centre. The flux vector over each edge of the triangular cell and the discrete form of the integral is
$$ \oint_{\gamma_{i}} \mathbb{F}(W)\cdot\textbf{n} d \sigma = \sum_{j\in N(i)}|\gamma_{ij}|\Phi(W_{ij},\textbf{n}_{ij}), $$
where,\ \	$ \displaystyle \Phi(W_{ij},\textbf{n}_{ij}) \simeq \frac{1}{|\gamma_{ij}|} \int_{\gamma_{ij}} \mathbb{F}(W)\cdot\textbf{n}_{ij} d \sigma,$ is the numerical flux computed at the interface between the cells $T_i$ and $T_j$.\\
As explained in the  \figurename{\ref{fig:mesh}}, $\gamma_{ij}$ is the edge surrounding the cell $T_i$ and $N(i)$ is the neighbouring triangles of the cell $T_i$. The intermediate solution $W_{ij}$ is reconstructed using the FVC scheme. The time discretization of (\ref{9}) is performed by a first order explicit Euler scheme. The time domain is divided into $N$ sub-intervals $[t_n , t_{n+1}]$ with time step $\Delta t  = t_{n+1} - t_n$ for $n = 0, 1,. . . . , N$. $W^n$ is the value of a generic function $W$ at time $t_n$. The fully-discrete formulation of the system (\ref{3}) is given by
\begin{equation}\label{10}
W^{n+1}_i=W^{n}_i-\frac{\Delta t}{|T_i|}\sum_{j\in N(i)}|\gamma_{ij}|\Phi(W^n_{ij},\textbf{n}_{ij})+\Delta t (S^n_i + Q^n_i ).
\end{equation}
\subsection{The FVC scheme}
In this subsection, we present a generalization of this scheme that was introduced in the preliminary works
\cite{benkhaldoun2015projection,ziggaf2020fvc}, the first was done in the Cartesian mesh and the second was done without taking into account the source term related to bathymetry and the balance produced by this term. For the corrector stage, we will use 2D finite volume formalism described in Section \ref{subsec31}. Finally,
the predictor stage and the final reformulation of the FVC scheme will be presented in this subsection.

\subsubsection{Construction of the projected speed model} 

Let us discretize the spatial domain $\Omega$  with cells $T_{i}$ as $\displaystyle \overline{\Omega}=\bigcup_{i=1}^{N_{ele}} T_{i}$ and  $\displaystyle \partial T_i=\bigcup_{j\in N(i)} \gamma_{ij}$,\\
with, \ 
$\partial T_i$ is the border of the cell $T_{i}$ and $N_{ele}$ is the total number of element.\\ 

Integrating the equations (\ref{1}) over the cell $T_{i}$, the basic equations of the finite volume method obtained using the divergence theorem are given by
\begin{subequations}\label{eq:sys1_eq}
	\begin{eqnarray}
	\displaystyle \frac{\partial}{\partial t} \int_{T_i} h \  \mathrm{d}V +  \int_{\partial T_i} hu{_\eta} \ \mathrm{d}\sigma &=&0 ,\label{eq:sys1_eq_1} \\[0.5em]
	\displaystyle \frac{\partial}{\partial t} \int_{T_i} hu \  \mathrm{d}V +  \int_{\partial T_i} \left\lbrace huu{_\eta}+ \frac{1}{2}gh^2 n_x  \right\rbrace \mathrm{d}\sigma &=& \displaystyle \int_{ T_i} -gh \partial_{x} Z \  dV + \int_{T_i} f_chv \ dV - \int_{T_i} \eta^2ghu\frac{|\textbf{u}|}{h^{4/3}} \ dV , \label{eq:sys1_eq_2} \\[0.5em]
	\displaystyle \frac{\partial}{\partial t} \int_{T_i} hv \ \mathrm{d}V +  \int_{\partial T_i} \left\lbrace hvu{_\eta}+ \frac{1}{2}gh^2 n_y  \right\rbrace \mathrm{d}\sigma &=& \displaystyle \int_{T_i} -gh \partial_{y}Z \  dV - \int_{T_i} f_chu \ dV - \int_{T_i} \eta^2ghv\frac{|\textbf{u}|}{h^{4/3}} \ dV , \label{eq:sys1_eq_3} 
	\end{eqnarray}
\end{subequations}\\
\begin{wrapfigure}[6]{r}{9cm}
	\vspace{-.9cm}
	\includegraphics[width=9.3cm,height=3.3cm]{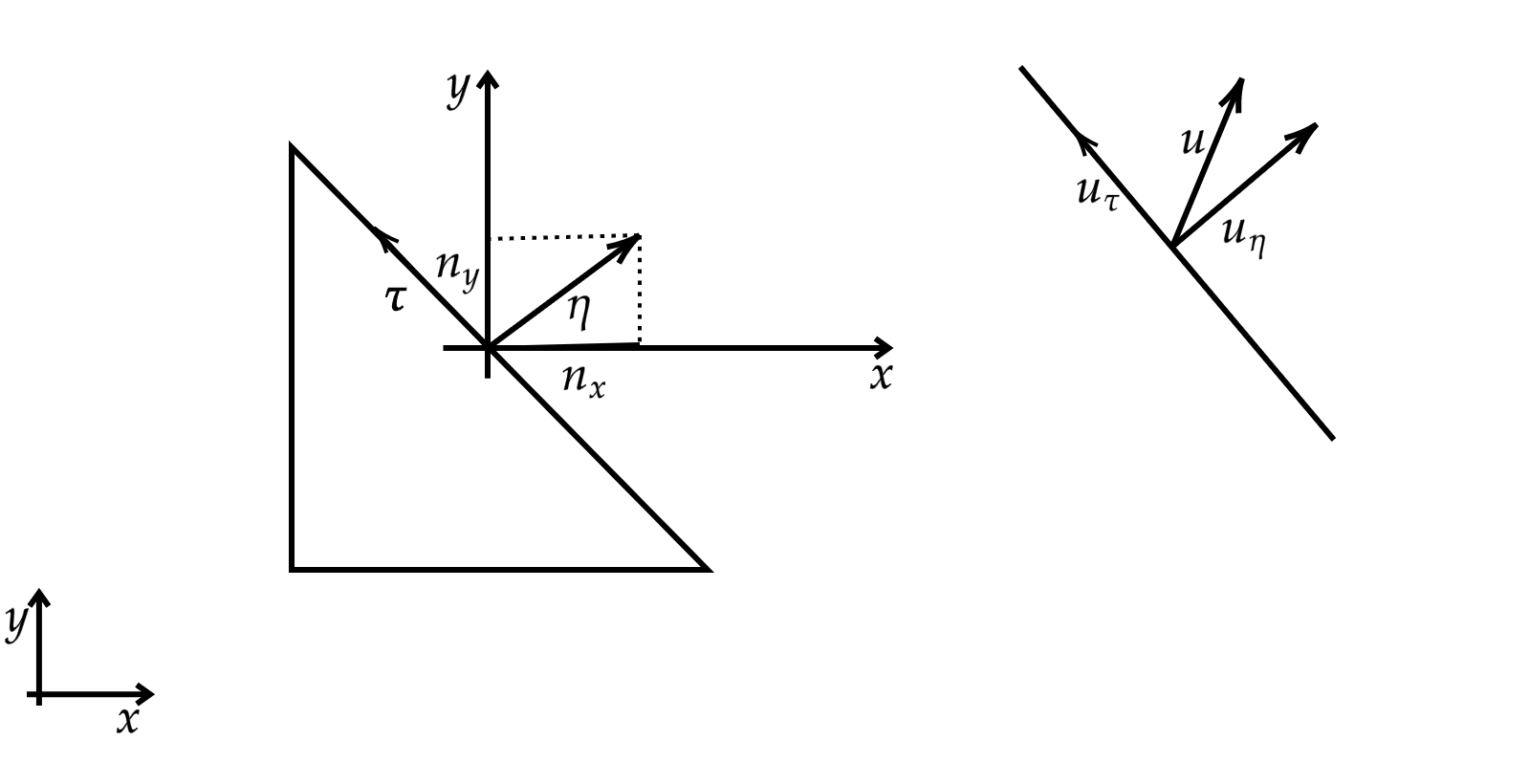}
	\caption{The projected velocity on the control volume.}
	\label{fig:project}
\end{wrapfigure}
where  $\eta = (n_x , n_y )^T$ the unit outward normal to the surface $|T_i|$ of
the cell $T_i$ and tangential $\tau = (-n_y , n_x )^T $ 
where the normal velocity $u_{\eta} = un_x + vn_y$ and the tangential velocity $u_{\tau} = vn_x -un_y$ (see \figurename{\ref{fig:project}}). In order to simplify the system (\ref{eq:sys1_eq}), we do the following operations\\[0.5em]
$\ref{eq:sys2_eq_2} \leftarrow n_x\ref{eq:sys1_eq_2} + n_y \ref{eq:sys1_eq_3}\ ,  \quad$ 
$\ref{eq:sys2_eq_3} \leftarrow n_x\ref{eq:sys1_eq_3} - n_y \ref{eq:sys1_eq_2}. $\ \ \
The outcome of these operations is\\[.em]

\begin{subequations}\label{eq:sys2_eq}
	\begin{eqnarray}
	\displaystyle \frac{\partial}{\partial t} \int_{T_i} h \  \mathrm{d}V +  \int_{\partial T_i} hu{_\eta} \ \mathrm{d}\sigma &=&0,\label{eq:sys2_eq_1} \\[0.5em]	\displaystyle \frac{\partial}{\partial t} \int_{T_i} hu_{\eta} \  \mathrm{d}V +  \int_{\partial T_i} \left\lbrace hu{_\eta}^2+ \frac{1}{2}gh^2  \right\rbrace \mathrm{d}\sigma &=& -\int_{ T_i} gh \nabla Z \cdot\textbf{n} \  dV + \int_{T_i} f_chu{_\tau}\ dV- \int_{T_i} \eta^2ghu{_\eta}\frac{|\textbf{u}|}{h^{4/3}} \ dV,\label{eq:sys2_eq_2}\\[0.5em]	\displaystyle \frac{\partial}{\partial t} \int_{T_i} hu_{\tau} \ \mathrm{d}V +  \int_{\partial T_i} hu_{\tau}u{_\eta} \  \mathrm{d}\sigma &=& - \int_{T_i} f_chu{_\eta}\ dV - \int_{T_i} \eta^2ghu{_\tau}\frac{|\textbf{u}|}{h^{4/3}} \ dV,\label{eq:sys2_eq_3}
	\end{eqnarray} 
\end{subequations}
which can be rewritten the system (\ref{eq:sys2_eq}) in a differential form as\\
\begin{equation}\label{13}
\left\{ \begin{array}{ll}
\displaystyle \frac{\partial h}{\partial t}   + \frac{\partial  hu{_\eta} }{\partial \eta}  =0,\\[1em] 
\displaystyle \frac{\partial hu{_\eta}}{\partial t}  +   \frac{\partial   }{\partial \eta}\left(  hu{_\eta}^2+ \frac{1}{2}gh^2  \right)   =-gh\partial_{\eta} Z+ f_chu_{\tau} - \eta^2ghu{_\eta}\frac{|\textbf{u}|}{h^{4/3}} ,\\[1em] 
\displaystyle \frac{\partial hu{_\tau}}{\partial t} +  \frac{\partial   }{\partial \eta}\left(  hu{_\eta}u{_\tau} \right) =-f_chu_{\eta} -\eta^2ghu{_\tau}\frac{|\textbf{u}|}{h^{4/3}}.
\end{array}
\right.
\end{equation}\\
The system (\ref{13}) can also be reformulated in the transport equation form as\\
\begin{equation}\label{14}
\displaystyle \frac{\partial \mathbf{U}}{\partial t}(t,X) + u_{\eta}(t,X) \frac{\partial\mathbf{U}}{\partial \eta}(t,X)= \mathbf{F}(\mathbf{U},Z,f_c) \quad  \quad  \forall \ \ X =(x,y)\in \Omega \subset \mathbb{R}^2,\ \ t>t_0,
\end{equation}\\
where,
$$\mathbf{U}= \left( \begin{array}{c} h\\[0.5em] hu_{\eta}\\[0.5em] hu_{\tau}\\ \end{array}\right), \quad
\textbf{F}(\mathbf{U},Z,f_c)= \left( \begin{array}{c} -h\partial_{\eta} (u_{\eta}) \\[0.5em] -gh\partial_{\eta} (h+Z)+f_chu_{\tau} - hu_{\eta}\partial_{\eta}(u_{\eta}) - \eta^2ghu{_\eta}\frac{|\textbf{u}|}{h^{4/3}}\\[0.5em] -f_chu_{\eta}-hu_{\tau}\partial_{\eta}(u_{\eta}) - \eta^2ghu{_\tau}\frac{|\textbf{u}|}{h^{4/3}}\\ \end{array}\right).$$ \\

The aim of using this technique on the local coordinates of the two-dimensional shallow water system (\ref{1}) in the control volume $T_{i}$, is  to reduce the  dimension to an one-dimensional system (\ref{14}) on each surface $|T_{i}|$ of this control volume.

\subsubsection{Flux construction}
we reconstruct the numerical flux $\Phi(W^n_{ij},\textbf{n}_{ij})$ using the method of characteristics. The fundamental idea of this method is to impose a regular grid at the new time level and to backtrack the flow trajectories to the previous time level, for more details see \cite{roe1986characteristic,seaid2001quasi}. At the previous time level, the quantities that are needed are evaluated by interpolation from their known values on a irregular grid, we'll see about that later.\\
\subsubsection{Method of characteristics}

The characteristic curves associated with the equation (\ref{14}) are solutions of the following Cauchy problem\\
\begin{equation}\label{15}
\left\{
\begin{array}{ll}
\displaystyle \frac{dX^c(t)}{d t} = u_{\eta}( t,X^c(t))\cdot\textbf{n} \quad  t \in [t_n, t_n +\alpha\Delta t],\ \ \alpha  \in ]0, \frac{t_{end}-t_n}{\Delta t} [,\\[0.8em] 

\displaystyle X^c(t_n + \alpha\Delta t)= X^*.
\end{array}
\right.
\end{equation}\\
\begin{wrapfigure}[16]{r}{7cm}
	\vspace{-.9cm}
	\includegraphics[width=7.5cm,height=4.9cm]{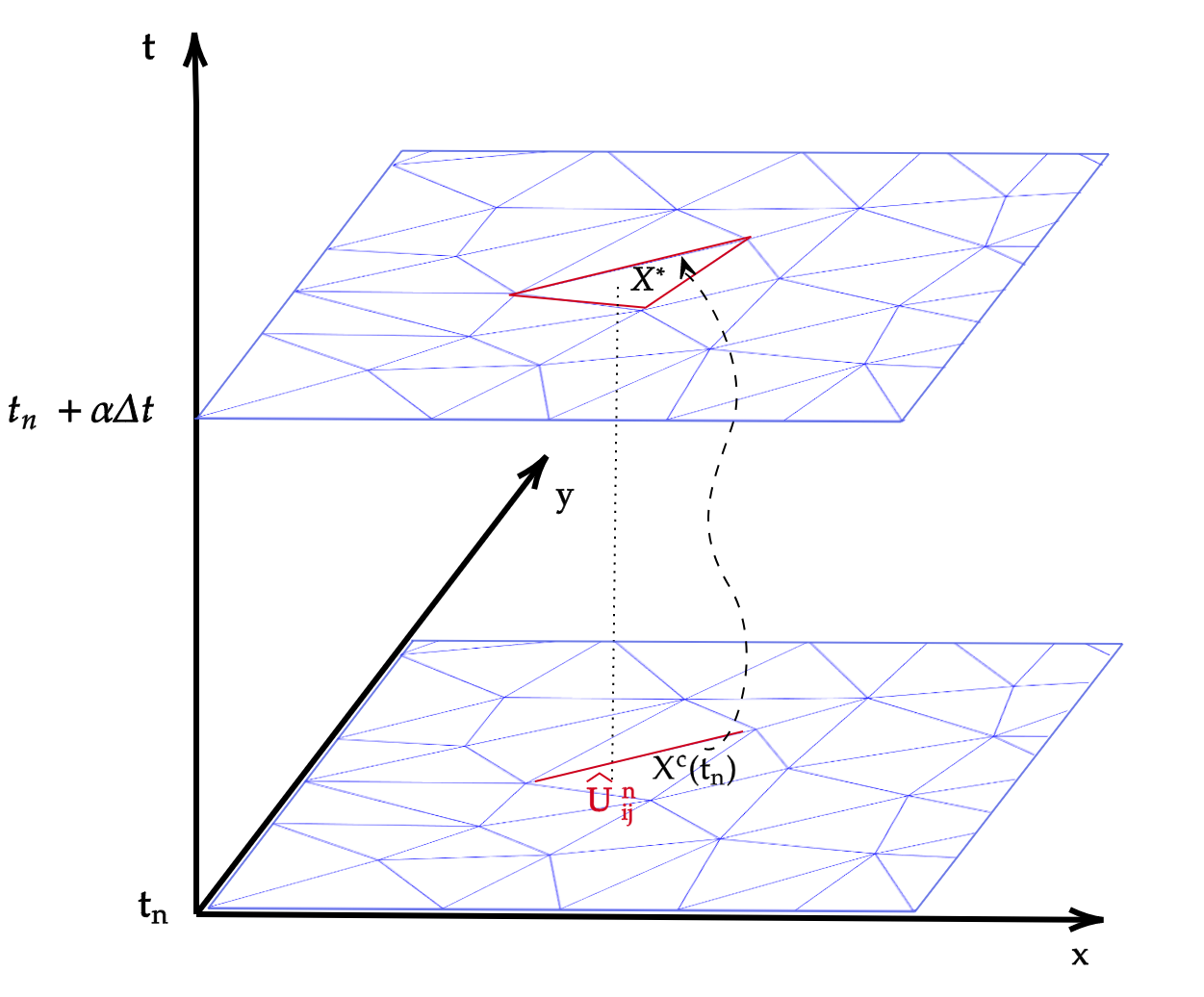}
	\caption{Illustration of the method of characteristics: An Eulerian gridpoint $X^c(t_n)$ is traced back in time to $X^*$ where the intermediate solution $\hat{\textbf{U}}^n_{ij}$ is interpolated.}
	\label{fig:meshtime}
\end{wrapfigure}
Note that $X^c(s)$ is the departure point at time $s$ of a particle that
will arrive at the interface $\gamma_{ij}$  in time $t_n + \alpha \Delta t$, see the \figurename{\ref{fig:meshtime}}. The method of characteristics does not follow the flow particles forward in time, as the Lagrangian schemes do, instead it traces backwards the position at time $t_n$ of particles that will reach the points of a fixed mesh
at time $t_n + \alpha\Delta t$. By doing so, the method avoids the grid distortion difficulties that the conventional Lagrangian schemes have. Hence, the solution of (\ref{15}) can be expressed in an integral form as\\
\begin{equation}\label{16}
X^c(t_n) = X^* - \int_{t_n}^{t_n + \alpha\Delta t} u_{\eta}(s,X^c(s))\cdot\textbf{n}\ d s,
\end{equation}
this integral can be calculated using an integral approximation method. In our approach, we used a Runge-Kutta 3 method to approximate the integral in (\ref{16}) which is accurate enough. In order to complete the reformulation of the algorithm used, the departure points must be calculated once the characteristic curves are known.
Therefore, the solution of the transport equation (\ref{14}) is given by\\
\begin{equation}\label{17}
\textbf{U}(t_n+ \alpha\Delta t,X^*) =\textbf{U}(t_n, X^c(t_n)) + \int_{t_n}^{t_n + \alpha\Delta t} \textbf{F}(	\textbf{U}(s,X^c(t_n)),Z,f_c)\ d s,
\end{equation}\\
where $	\textbf{U}(t_n,X^c(t_n))$  is the solutions at the characteristic feet computed by the local least squares interpolation method. In other cases, the integral of the equation (\ref{17}) can be calculated using a first-order approximation based on the rectangle method, the vector $\textbf{U}^n$ is reconstructed at the interfaces using 
\begin{equation}\label{18}
\displaystyle	\textbf{U}^n_{ij} = \hat{\textbf{U}}^n_{ij} + \alpha \Delta t \textbf{F}( \hat{\textbf{U}}^n_{ij},Z,f_c),
\end{equation}
where,\\
\begin{equation}\label{19}
  \hat{\textbf{U}}^n_{ij} = \textbf{U}(t_n,X^c(t_n)) = \sum_{k \in V(c)} \alpha_k(c)\textbf{U}(X^k),
\end{equation}

\begin{wrapfigure}[7]{r}{7cm}
	\vspace{-.9cm}
	\includegraphics[width=7.5cm,height=3.5cm]{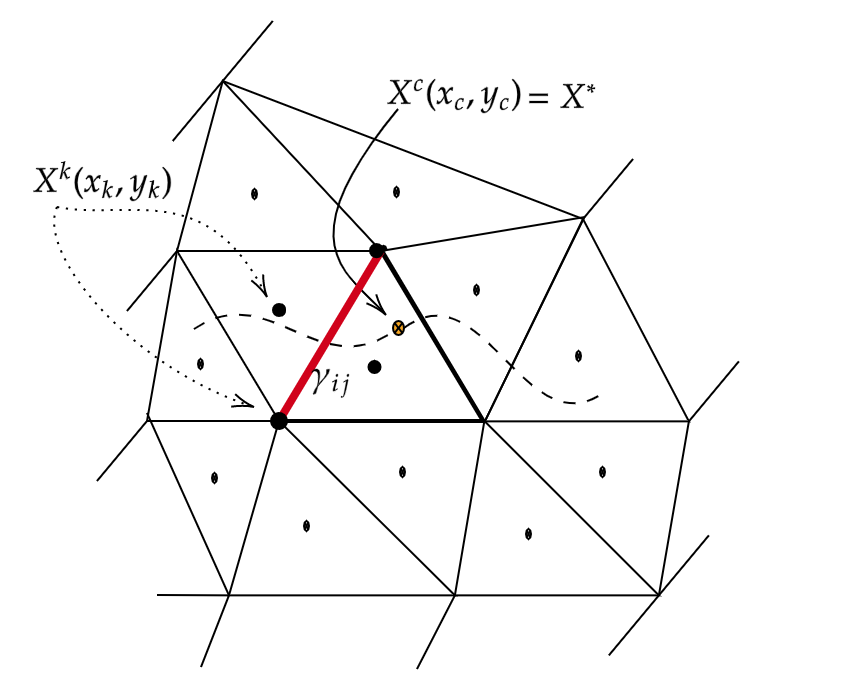}
	\caption{LSM Illustration.}
	\label{fig:LSM}
\end{wrapfigure}
where, $V(c) := \{$ the nodes and cell centres around the edge $ \gamma_{ij}\} $, see the \figurename{\ref{fig:LSM}} and $ \alpha_k(c) $ is weights coming from the least squares method (LSM). It  can be written
\begin{equation}\label{20}
\alpha_k(c) = \frac{1 + \lambda \cdot(X^k -X^c)}{Card(V(c)) + \lambda\cdot R} ,
\end{equation}
 such as, $\lambda=(\lambda_x,\lambda_y)$ and  $R=(R_x,R_y)$. The weights parameters are given by formulas (see the Appendix).\\[2em]

\begin{wrapfigure}[13]{r}{6.1cm}
	\vspace{-.5cm}
	\includegraphics[width=6cm,height=5.cm]{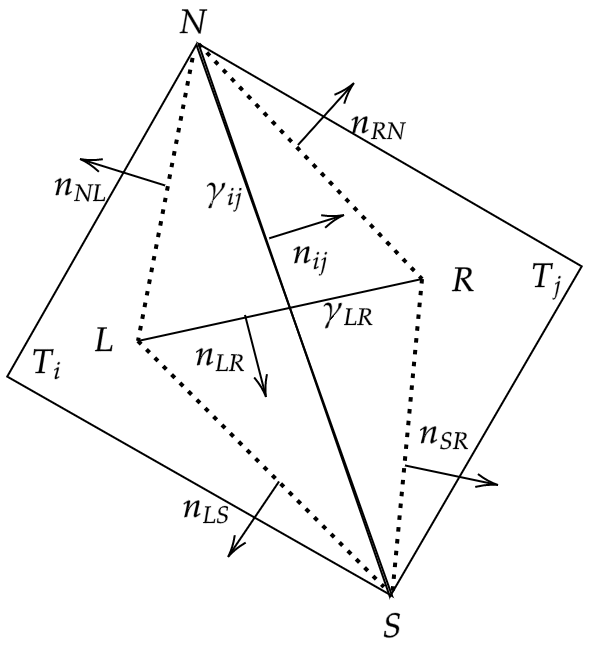}
	\caption{Diamond cell in 2D}
	\label{fig:Diamond}
\end{wrapfigure}
To approximate $\textbf{F}(\textbf{U},Z,f_c) $, (i.e. $ \partial_{\eta}(u_{\eta})$,  $ \partial_{\eta}(h + Z),...$ ) we will need to approximate these derivatives in the interfaces,  for that we use the diamond cell \figurename{\ref{fig:Diamond}}. This cell is constructed by connection of centres of gravity $(L, R)$ of cells $T_i$, $T_j$ which share the interface $\gamma_{ij}$ and its endpoints $S$, $N$. We obtain the co-volume $SRNL$ by this construction.
One can assume that the gradient is constant on the co-volume $SRNL$. According
to Green-Gauss theorem the approximation leads to
\begin{equation}\label{21}
\nabla u_{ij} =\frac{1}{2\mu_{_{SRNL} }} \left\lbrace \frac{}{}(u_{_{S}}-u_{_{N}}) \vec{n_{_{LR}}}|\gamma_{_{LR}}| + (u_{_{R}}-u_{_{L}}) \vec{n_{ij}}|\gamma_{ij}|  \right\rbrace ,
\end{equation}
where $u_{_{N}}$, $u_{_{S}}$, $u_{_{R}}$, and $u_{_{L}}$  represent respectively the values of the quantity $u$ in the point $N$, $S$, $R$ and $L$.
$\vec{n_{_{LR}}}$ is a unit normal vector of the co-volume face $\gamma_{_{LR}}$ and $|\gamma_{_{LR}}|$ is its length. The others co-volume interfaces and their normal vectors are labeled analogically. $\mu_{_{SRNL} }$ is the area of the co-volume $SRNL$.\\

After the discretization of the source terms (see subsection  \ref{subsec325}), the district equation system (\ref{18}) leads to the following predictor step\\
\mathleft

\textbullet) \textit{Predictor stage}\\
\begin{equation}\label{22}
\hspace{0.6cm}\left| 
 \begin{aligned}
\displaystyle h^n_{ij}  &= \  \hat{h}^n_{ij} - \alpha \Delta t \hat{h}^n_{ij} \nabla  \hat{(u_{\eta} })^n_{ij} ,\\
\displaystyle (hu_{\eta} )^n_{ij} &= \hat{(hu_{\eta} })^n_{ij} - \alpha  \Delta t \left\lbrace  g  \hat{h}^n_{ij} \nabla(\hat{h}^n_{ij} + Z_{ij}) -  f_c \hat{(hu_{\tau} })^n_{ij} + \hat{(hu_{\eta} })^n_{ij}\nabla  \hat{(u_{\eta} })^n_{ij} + \eta^2g \hat{(hu_{\eta} })^n_{ij} \frac{|(\hat{\textbf{u}})^n_{ij}|}{(\hat{h}^{3/4})^n_{ij}}\right\rbrace, \\
\displaystyle (hu_{\tau} )^n_{ij} &= \hat{(hu_{\tau} })^n_{ij} - \alpha \Delta t\left\lbrace  f_c\hat{(hu_{\eta} )}^n_{ij} + \hat{(hu_{\tau} )}^n_{ij}\nabla  \hat{(u_{\eta} })^n_{ij} + \eta^2g \hat{(hu_{\tau} })^n_{ij} \frac{|(\hat{\textbf{u}})^n_{ij}|}{(\hat{h}^{3/4})^n_{ij}}\right\rbrace .
\end{aligned}
\right.
\end{equation}

Once these projected states are calculated, the quantity $W_{ij}$
will be  calculated using the following transformations\\ 

$hu^n_{ij} = (hu_{\eta} )^n_{ij}n_x -(hu_{\tau} )^n_{ij}n_y, \qand hv^n_{ij} = (hu_{\tau} )^n_{ij}n_x + (hu_{\eta} )^n_{ij}n_y $.\\

\textbullet) \textit{Corrector stage}\\ 
\begin{equation}\label{23}
\hspace{0.6cm}\left| 
\ \ \begin{aligned}
\displaystyle W^n_{ij} \  \   \ \  \ \  = \ \ \ & (
h^n_{ij} \ \ \  hu^n_{ij} \ \ \  hv^n_{ij}
)^T
,\\[.8em]
\displaystyle \Phi(W^n_{ij},\textbf{n}_{ij}) &= \mathbb{F}(W^n_{ij}) \cdot \textbf{n}_{ij}, \\[.8em]
\displaystyle W^{n+1}_i \ \  = \ \ \ & W^{n}_i-\frac{\Delta t}{|T_i|}\sum_{j\in N(i)}|\gamma_{ij}|\Phi(W^n_{ij},\textbf{n}_{ij}) +\Delta t S^n_i +\Delta t Q^n_i .
\end{aligned}
\right.
\end{equation}
\mathcenter
\subsubsection{Well-balanced FVC scheme: the discretization of the bathymetry source term} \label{subsec325}
 In order to be able to calculate realistic flows we now consider the case $\grad{Z} \neq 0_{\mathbb{R}^2} $ and introduce a numerical discretization of the source terms. As discussed in paragraph \ref{subsec211}, the treatment of source terms related to bathymetry in the shallow water system poses a challenge in many numerical methods. In our scheme, the approximation of the source term $S^n_i$ is reconstructed in such a way that the C-property \cite{bermudez1994upwind} is satisfied, i.e. to maintain a discrete local balance of the continuous stationary state in still water.\\
\begin{equation}\label{24}
\left.  \begin{aligned}
\displaystyle  h^n_i+Z_i=& h^n_j +Z_j = H := cst\\ \textbf{u}_i^n + \textbf{u}_j^n =& 0_{\mathbb{R}^2},\ \  \forall\  T_i, T_j \in \overline{\Omega}\\
\end{aligned}
\right\rbrace \Rightarrow h_i^{n+1} + Z_i = H, \qand \textbf{u}_i^{n+1}= 0_{\mathbb{R}^2}.
\end{equation}\\

\textbullet) \textit{The hydrostatic balance}\\
 
$\displaystyle
 \grad({\frac{1}{2}gh^2}) = -gh\grad{Z},
$\\
we prove from the hydrostatic balance that the model of the projected speed preserves the stationary state of the lake at rest

\begin{wrapfigure}[10]{r}{6.1cm}
	\includegraphics[width=6.5cm,height=4.8cm]{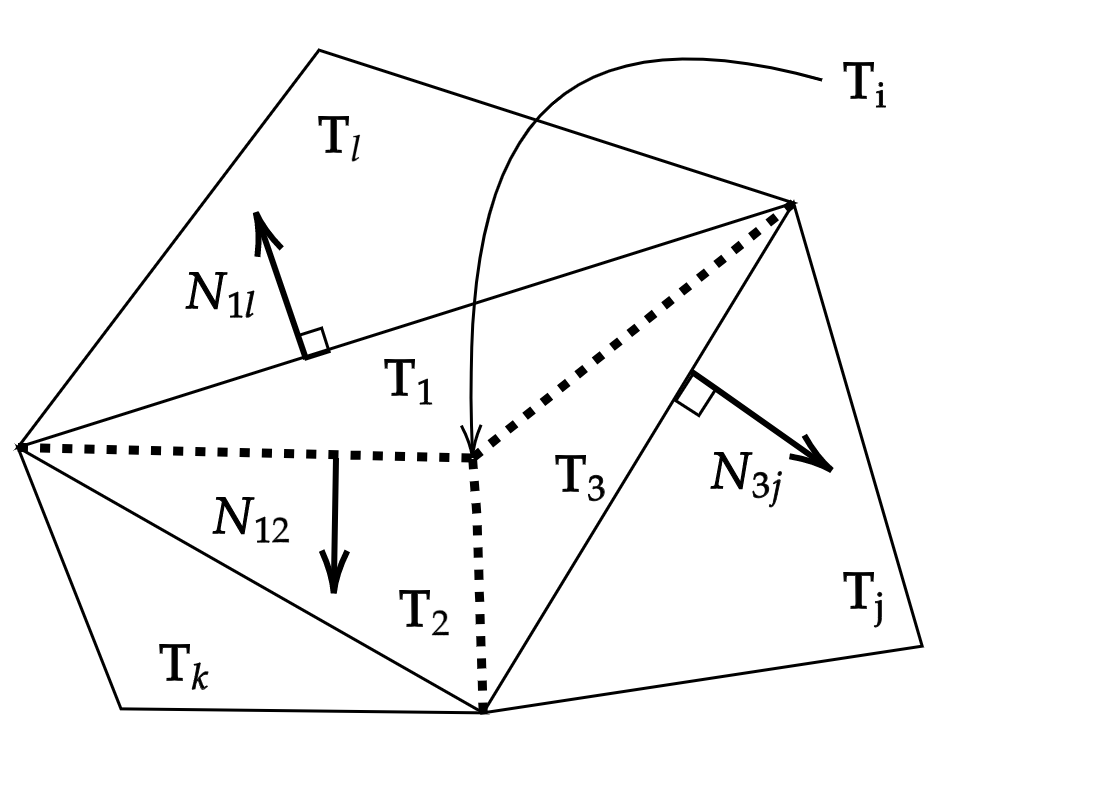}
	\caption{Sub-triangles used in the discretization of source terms}
     \label{fig:sourceterme}
\end{wrapfigure}

\textbullet) \textit{The projected speed model}\\
\begin{equation}\label{25}
\partial_{t}\left( \begin{array}{c} h\\[0.5em] 0\\[0.5em] 0\\ \end{array}\right)
+ 0 \times \partial_{\eta} \left( \begin{array}{c} h\\[0.5em] 0\\[0.5em] 0\\ \end{array}\right)
= \left( \begin{array}{c} 0\\[0.5em] -g\partial_{\eta}(h +Z)\\[0.5em] 0\\ \end{array}\right),
\end{equation}\\
$ 
\displaystyle \partial_{t} h =0, \qand \partial_{\eta}(h+Z) = 0 \Longrightarrow h(x,y,t) + Z(x,y) = cst\  \ \forall \ x,\  y,\ t .
$\\[.5em]
This result ensures the equilibrium property corresponding to the lake at rest, and therefore it is consistent with the continuous form of the system's equilibrium with bathymetry source term.\\

\begin{equation}\label{26}
\frac{\Delta t}{|T_i|}\sum_{j\in N(i)}|\gamma_{ij}|\Phi(W^n_{ij},\textbf{n}_{ij}) = \frac{\Delta t}{|T_i|} \int_{T_i} S dV,
\end{equation}\\
which is equivalent to\\
\begin{equation}\label{27}
\displaystyle
\left( \begin{array}{c} 0\\[0.5em] \sum_{j\in N(i)} \frac{1}{2}g(h_{ij})^2(n_{ij})_x|\gamma_{ij}| \\[0.5em] \sum_{j\in N(i)} \frac{1}{2}g(h_{ij})^2(n_{ij})_y|\gamma_{ij}|\\ \end{array}\right)
= 
\left( \begin{array}{c} 0\\[0.5em] -g\int_{T_i}h\partial_{x}Z dV \\[0.5em] -g\int_{T_i}h\partial_{y}Z dV\\ \end{array}\right).
\end{equation}\\
To approximate the source terms, we proceed as follows. First, we decompose the triangle $T_i$ into three sub-triangles, as depicted in \figurename
{\ref{fig:sourceterme}}.\\
where $ N_{xij}= (n_{ij})_x |\gamma_{ij}|,  \qand N_{yij}= (n_{ij})_y |\gamma_{ij}| $. Then, the source term is approximated as

\begin{equation}\label{28}
\int_{T_i} h \partial_{x} Z dV = h_1\int_{T_1} \partial_{x} Z dV +h_2\int_{T_2}  \partial_{x} Z dV +h_3\int_{T_3}\partial_{x} Z dV,\\
\end{equation}
with $h_1$, $h_2$ and $h_3$ are the average values of $h$ over $T_1$, $T_2$ and $T_3$ respectively. 

\begin{equation}\label{29}
\begin{aligned}
h_1\int_{T_1} \partial_{x} Z dV&=\sum_{j\in N(1)} \int_{\gamma_{1j}} Z n_x dV\\
&= h_1\sum_{j\in N(1)} \frac{Z_1+Z_j}{2}N_{x1j}\\
&= \frac{h_1}{2}\left\lbrace (Z_1 +Z_l)N_{x1l} + (Z_1 +Z_2)N_{x12}  +(Z_1 +Z_3)N_{x13} \right\rbrace.\\
  \end{aligned}
\end{equation}
The same applies to the y-direction.\\ Again the stationary flow condition $\displaystyle h_1 + Z_1 = h_j+ Z_j = H = cst$, $\forall j\in N(1)$ $\Rightarrow$
$\displaystyle h_1 +h_j +Z_1 +Z_j = 2H$ and $ H - \frac{h_1 + h_j }{2} = \frac{Z_1 +Z_j}{2}$. Thus, (\ref{29}) gives \\

$\displaystyle \int_{T_1} h\partial_{x} Z dV = h_1\sum_{j\in N(1)}\left( H - \frac{h_1 +h_j}{2} \right) N_{x1j}  \underset{\displaystyle \overbrace{\sum_{j\in N(1)} N_{x1j}=0}}{=} -\frac{h_1}{2}\sum_{j\in N(1)} h_j N_{x1j} $ . \\
Finally, \\

$\displaystyle \int_{T_1} \partial_{x} Z dV =-\frac{h_1}{2}(h_lN_{x1l}+h_2N_{x12}+h_3N_{x13})$ .\\ 


A similar procedure leads to the following approximations of the other terms in (\ref{28})

$\displaystyle \int_{T_2} \partial_{x} Z dV =-\frac{h_2}{2}(h_kN_{x2k}+h_1N_{x21}+h_3N_{x23})$,\\
 
$\displaystyle \int_{T_3} \partial_{x} Z dV =-\frac{h_3}{2}(h_jN_{x3j}+h_1N_{x31}+h_2N_{x32})$.\\[0.5em]
Notice that $h_l , h_k$ and $h_j$ are the average values of $h$, respectively, on
the triangle $T_l, T_k$ and $T_j$. Summing up, using the fact that  ($N_{xij} = -N_{xji}$) so,  the discretization (\ref{28}) gives\\
\begin{equation}\label{30}
\displaystyle  \int_{T_i} h\partial_{x} Z dV = -\frac{1}{2} \left(  h_1h_lN_{x1l}+ h_2h_kN_{x2k}+ h_3h_jN_{x3j}\right).
\end{equation}\\
For this reconstruction, the source terms in (\ref{27}) result in\\

\begin{equation}\label{31}
\begin{aligned}
\sum_{j\in N(i)}(h^n_{ij})^2 N_{xij} =& h_1h_lN_{x1l}+ h_2h_kN_{x2k}+ h_3h_jN_{x3j}, \\[.5em]
\sum_{j\in N(i)}(h^n_{ij})^2 N_{yij} =& h_1h_lN_{y1l}+ h_2h_kN_{y2k}+ h_3h_jN_{y3j} .
\end{aligned}
\end{equation}

If you have noticed, we will need $h_1, h_2$ and $h_3$ to be able to calculate the values of the integrals in equation (\ref{30}) but the system. (\ref{31}) has two equations for the three unknowns. To complete the system, we add the natural conservation equation, $\displaystyle h_1 + h_2 + h_3 = 3h_i $. The following system gives us the values we need  
\begin{equation}\label{32}
\left( \begin{array}{c}
h_1\\ [.5em] 
h_2\\ [.5em] 
h_3
\end{array}\right) = \left( 
\begin{array}{ccc}
1	& 1 & 1 \\[.5em] 
h_lN_{x1l}	&h_kN_{x2k}  & h_jN_{x3j}  \\ [.5em] 
h_lN_{y1l}	&h_kN_{y2k}  & h_jN_{y3j} 
\end{array}\right)^{-1}\cdot\left( \begin{array}{c}
3h_i\\ [.5em] 
\sum_{j\in N(i)}(h^n_{ij})^2 N_{xij}\\ [.5em] 
\sum_{j\in N(i)}(h^n_{ij})^2 N_{yij}
\end{array}\right).
\end{equation}

Analogously, the bottom values $Z_j$ , $j = 1, 2, 3 $ are reconstructed in
each sub-triangle of $T_i$ as\\

$\displaystyle\left\lbrace  \begin{array}{c}
Z_1	 = Z_i +h_i^n - h_1^n,\\[.5em] 
Z_2	 = Z_i +h_i^n - h_2^n,\\[.5em] 
Z_3	 = Z_i +h_i^n - h_3^n.
\end{array}\right.$\\ 

Finally the source terms in (\ref{30}) are approximated as\\

\begin{equation}\label{33}
\int_{T_i} h\partial_{x} Z dV = h_1 \sum_{m\in N(1)} \frac{Z_1 - Z_m}{2}N_{x1m} +h_2 \sum_{m\in N(2)} \frac{Z_2 - Z_m}{2}N_{x2m}  +h_3 \sum_{m\in N(3)} \frac{Z_3 - Z_m}{2}N_{x3m},
\end{equation}

with a similar equation for the other source terms in the y-direction.\\


\textbullet) \textit{ Computation of the solution}\\

 Finally, we write the formally well-balanced FVC scheme after calculation of the  interface values (\ref{23}) and the bathymetry source term approximation (\ref{33}) as 
\begin{equation}\label{34}
 W^{n+1}_i = W^{n}_i-\frac{\Delta t}{|T_i|}\sum_{j\in N(i)}|\gamma_{ij}|\Phi(W^n_{ij},\textbf{n}_{ij}) +\Delta t S^n_i.
\end{equation}
\subsubsection{Semi implied  treatment of  friction term  source }

To avoid stability problems related to the bottom friction source term, a fractional semi implied treatment for this term is proposed. The idea is to evaluate the momentum in the system (\ref{3}) by decomposing it into two equations

\begin{equation}\label{36}
\left\{
\begin{array}{ll}
\displaystyle \frac{\partial h\textbf{u}}{\partial t} = -\eta^2gh^{-1/3}|\textbf{u}|\textbf{u} \\[0.8em] 

\displaystyle \frac{\partial h\textbf{u}}{\partial t} + \Phi_{h\textbf{u}}(W)= -gh \grad{Z},
\end{array}
\right.
\end{equation}
where $\Phi_{h\textbf{u}}(W)$ represents the convection terms corresponding to the equations of the momentum.
In a first step, a linearized semi implied method is used to integrate the first equation of the system (\ref{36})
\begin{equation}\label{37}
	\frac{(\tilde{h\textbf{u}})_i - (h\textbf{u})^n_i}{\Delta t} = - \eta^2g(\tilde{h\textbf{u}})_i|\textbf{u}^n_i|(h^n_i)^{-4/3}.
\end{equation}
In the second step, the value $(\tilde{h\textbf{u}})_i$ is taken as the initial condition for solving the second  equation of (\ref{36}).
%
\section{Numerical tests} \label{sec4}

It is clear from (\ref{23}) that the scheme is conservative  and can compute the numerical flux corresponding to the physical solutions of water flow without relying on Riemann problem solvers. The CFL condition for the explicit scheme (\ref{22}) can be written\\
$\displaystyle \Delta t \leqslant \min \left\lbrace  \frac{|T_i|}{|\gamma_{ij}|(\textbf{u}\cdot\textbf{n} + \sqrt{gh})} , \frac{|T_i|}{|\gamma_{ij}|(\textbf{u}\cdot\textbf{n} + \sqrt{gh})\sqrt{2\alpha}}  \right\rbrace $.

A fixed  CFL = 0.9  is used  and $\alpha = 2 $ if not specified in the test. The used computer is an Intel Core i7-8565U CPU @ 1.80GHz $\times$ 8, with  15  GB RAM.
 
In order to validate our FVC scheme on unstructured meshes to simulate shallow water flows, therefore we present some test cases that are proposed by several authors to validate their model and their numerical approach. The accuracy is demonstrated by comparing numerical solutions produced by the FVC scheme with analytical solutions, especially in tests  \ref{Accuracy}, \ref{Dambrak} and \ref{irrebed}.  
 To reproduce
the calculation results reported in the literature, the source term of the bed is always taken into account. The C-property produced by this term has also been treated in tests \ref{Cprop} and.
The Coriolis effect was taken into account in test \ref{Circular} and this test’s results are in  good agreement with those presented in the literature.

\subsection{ Accuracy test example}\label{Accuracy}
\begin{wrapfigure}[18]{r}{0cm}
    \includegraphics[width=.5\linewidth,height=0.25\textheight]{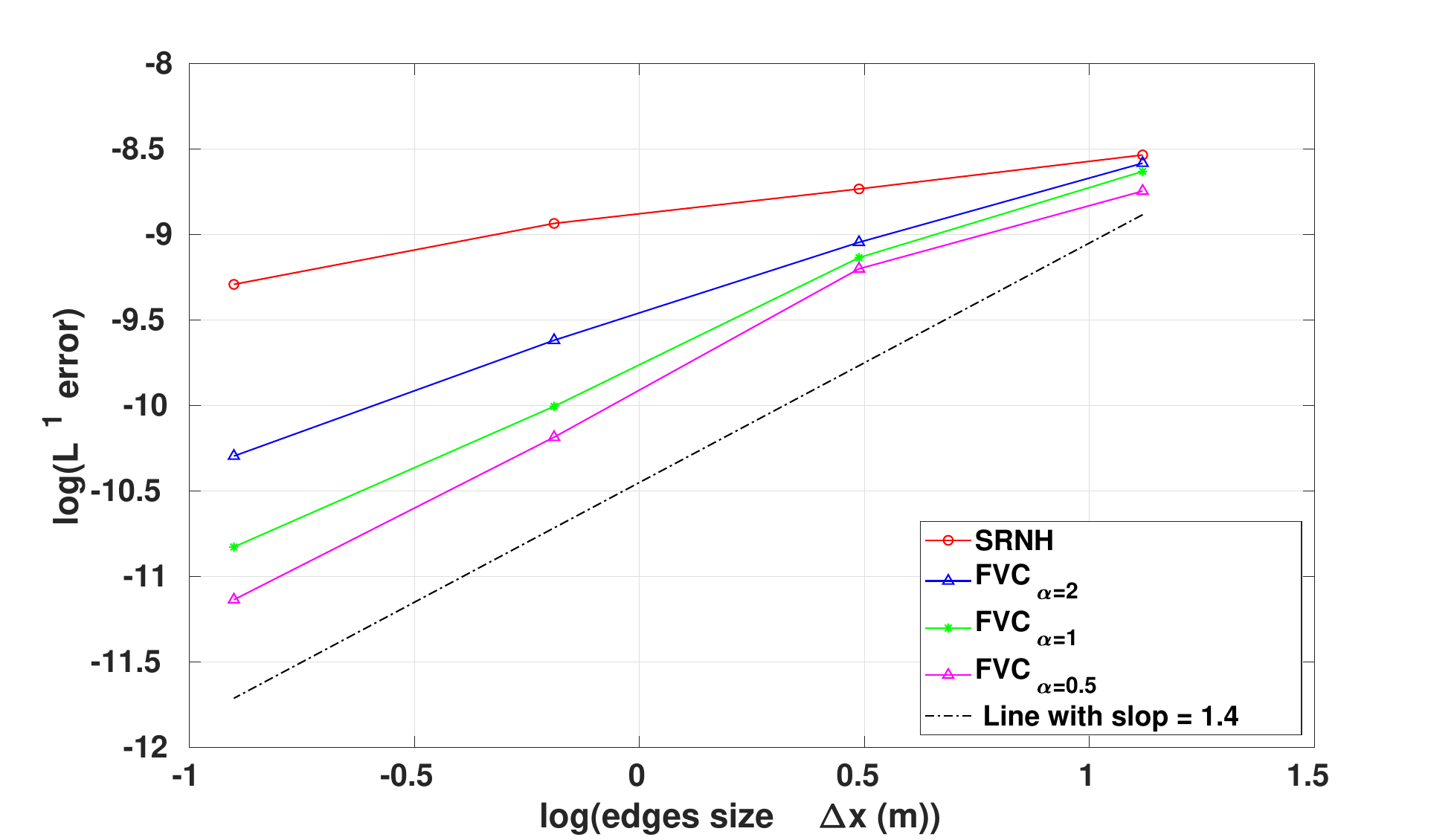}
	\caption{Convergence order in L$^1$ error of a water height.}
    \label{fig:convergence}
\end{wrapfigure}

We test our approach on a problem where the exact solution is known \cite{fjordholm2011vorticity}. It can be readily checked that.
\mathleft 
\begin{equation*}\label{35}
\begin{aligned}
 h(t,x,y) &= 1 -\frac{a^2}{4bg} \exp(-2b(\bar{x}^2 + \bar{y}^2)), \\[.2em]
 u(t,x,y) &= Mcos(\theta) + a\bar{y}\exp(-b(\bar{x}^2 + \bar{y}^2)),\quad v(t,x,y) = Msin(\theta) - a\bar{x}\exp(-b(\bar{x}^2 + \bar{y}^2)),
\end{aligned}
\end{equation*}
where, \ \ 
$ \displaystyle \bar{x} = x - x_0 -Mt\cos(\theta)$ and $ \displaystyle \bar{y} = y - y_0 -Mt\sin(\theta)$,\\
 
gives a smooth solution of the shallow water system (\ref{1}) without the source terms ( i.e. $\grad{Z} = 0, \qand f_c = 0$  ) for any
choice of constants, $M, a, b, x_0, y_0 $ and $\theta$. Initial and boundary condition are set according to the exact solution.\\ We let $M = \frac{1}{2},\  g = 1,\  a = 0.04,\  b = 0.02$, and $(x_0, y_0) = (-20,-10)$. To test the scheme ability
to resolve flows that are not aligned with the computational mesh, we let $\theta = \frac{\pi}{6}$.
We compute for $(x,y) \in \Omega = [-50, 50] \times [-50, 50]$ up to time $t = 100 $ s to compare the accuracy of our FVC scheme to the SRNH scheme originally introduced in \cite{monthe1999positivity,roe1981approximate}. We also propose to see the effect of the choice of the parameter $\alpha$. In \cite{benkhaldoun2010simple} where the authors defined for the first time the FVC scheme in its one-dimensional formulation, were able to show that the parameter $\alpha$ controls the accuracy of the FVC scheme, for more rigorous details we invite you to see the Lemma 3.2 and its proof. Our study reaches the same conclusion on alpha as in \figurename{\ref{fig:convergence}} we show the convergence order for three different choices of $\alpha$.  
\mathcenter 
\begin{table}[H]
\caption{Relative $L^1$ errors and CPU times obtained for the accuracy test example at time t = 100 using the SRNH and FVC schemes.}\label{Tab1} 
\centering
\begin{tabular}{lccccccccl}\toprule
Schemes	& \multicolumn{3}{c}{SRNH} & \multicolumn{6}{c}{FVC$_{\alpha=1}$}
	\\\cmidrule(lr){2-5}\cmidrule(lr){7-10}
	& \multicolumn{3}{c}{$L^1$ error } & CPU time (s) &  &\multicolumn{3}{c}{$L^1$ error }&CPU time (s)
    \\\cmidrule(lr){2-4}\cmidrule(lr){7-9}
	$\#$ Cells & $h$       & $hu$      & $hv$     & \      &\ &$h$       & $hu$     & $hv$    &  \\\midrule
	2648       & 2.042E-04 &5.344E-03  & 9.546E-04& 16.71  &\ &1.785E-04 & 4.379E-03&8.763E-03& 11.91\\
	10362      & 1.627E-04 &7.626E-03 & 7.626E-03& 31.30  &\ &1.078E-04 & 1.534E-03&3.475E-03& 21.49\\
	40690      & 1.317E-04 &2.716E-03 & 5.450E-03& 184.81 &\ &4.518E-05 & 5.430E-04&1.249E-03& 101.46\\
    161316     & 9.363E-05 &1.628E-03  & 3.441E-03& 1887.89&\ &1.986E-05 & 2.083E-04&4.868E-04& 1003.26\\
     640138    & 8.526E-05 &9.871E-04  & 1.067E-03& 21031.82&\ &1.066E-05 & 1.771E-04&3.278E-04& 10524.03\\\bottomrule
\end{tabular}
\end{table}

The results  in Table \ref{Tab1} show that increasing the number of the cells in the computational domain lead to a decrease in the $L^1$ error for the water height $h$, and the discharges $hu$ and $hv$ in all schemes. Faster decay of the error is observed in the FVC scheme than in the SRNH scheme which is only natural a simple inspection of Table \ref{Tab1} also reveals that for meshes with a low number of the cells, the measured computation time is comparable for the SRNH scheme and FVC scheme. 
However, for meshes with a somewhat large number of cells, the FVC method is the most efficient. For example, a mesh of $161316$ cells, the FVC scheme is about two faster than the SRNH scheme this is due to our FVC approach is not based on the calculation of the Jacobian matrix of the system, this matrix intervenes in many Q-scheme type approximation schemes, it is responsible for the slowness of this kind of scheme. Note that the SRNH scheme requires a solver for the Riemann problem at each time step to reconstruct the numerical flux, which is completely avoided in our FVC scheme. We can see that both schemes could reach the designed order of accuracy. In \figurename{\ref{fig:convergence}} we have plotted the $\log$ of $L^1$ error calculated in Table \ref{Tab1} and the other values calculated for $\alpha = 0.5$ and then $\alpha = 2$ against the $\log$ of the maximum value of the mesh edges.
We find that the $ L ^ 1 $ errors of the FVC scheme lie on a slope line $1.4$, indicating that the accuracy order of the scheme is about $1.4$ for $ \alpha = 0.5$.

\subsection {Dam-break problem}\label{Dambrak}
Flood flows produced by the dam break, segments of dykes, or other structures are torrential in nature with the presence of a discontinuous front propagating downstream and a rarefaction wave propagating upstream. The characteristics of these flows such as velocity, water level and time of arrival of floods must be determined in advance in order to manage floods and reduce their impact on the environment and economic infrastructure.
In order to test our approach for problems related to dam break, we carried out a series of test cases proposed in the literature ( see, e.g. \cite{benkhaldoun2010simple,dafermos2005hyperbolic,stoker1948formation} ) etc.\\
The proposed approach is not based on a Riemann solver technique which is very appropriate for a hyperbolic problem whose solution is often represented by a discontinuous front. Thus it will be interesting to examine the accuracy of our approach by simulating torrential flows with the presence of a discontinuity in the velocity profile and the free surface.\\[-0.2cm]

\textbullet) \textit{ Description of the problem}\\[-0.2cm]

We consider a rectangular channel with a flat bottom, $Z(x,y) = 0$ and no friction, i.e. there is no source terms, the problem is purely hyperbolic. The channel is $1.6\  m$ long and $0.1\  m$ wide (we assume a nondimensionalization problem), the initial conditions are given by
\mathleft
\begin{equation*}
h(0,x,y) = \left\{ \begin{array}{ll} 
h_l  \quad \quad \quad if \ x \leqslant x_m,\\[.5em]
h_r  \quad \quad \quad if  x > x_m,
\end{array}\right.  \quad    0 \leqslant y \leqslant 0.1
\end{equation*}
$ u(0,x,y) = v(0,x,y) = 0\  m/s.$ \  \ A dam is placed in the middle of the channel i.e. $x_m = 0.8 m$.\\
\mathcenter
This corresponds to a homogeneous Riemann problem. Initially, the water is at rest, the height $h_l$ remains $1.0$ m for all simulations. The downstream height $h_r$ takes on different values $0.5$ m, $0.1$ m, and $0.025$ m. The nature of the torrential flow due to the dam break depends essentially on the ratio $h_r / h_l$.\\
A numerical instability is likely to occur for small values of the ratio $h_r/h_l$. At $t = 0$, it is assumed that the dam is abruptly removed causing a shock wave with the presence of a discontinuous front of the water surface propagating downstream.\\
The channel is assumed to be closed on all four sides and the "slip" conditions is imposed on all walls. The computational domain is discretised by a mesh of $41776$ triangles with an average size of $0.003 m^2$. 
We will compare the water height and velocity obtained by our approach with the analytical solution which is calculated using the Stoker method \cite{stoker2011water}.

The evolution of the water surface profile is used to examine the
behaviour of the FVC scheme in capturing the discontinuous shock front. Henderson \cite{henderson1966macmillan} notes that when the ratio $h_r/h_l$ is greater than 0.138, the flow is sub-critical in the whole of the channel. When the ratio $h_r/h_l$ is smaller than 0.138, the flow is supercritical downstream and sub-critical upstream of the dam.\\
For very small values of $h_r/h_l$, the upstream flow regime becomes strongly supercritical, and it may be difficult to capture such a shock wave numerically.\\[-0.2cm]

\textbullet) \textit{ Results and discussion}\\[-0.2cm]

\begin{figure}[H] 
	\begin{subfigure}[b]{0.5\linewidth}
		\centering
		\includegraphics[width=.99\textwidth,height=.2\textheight]{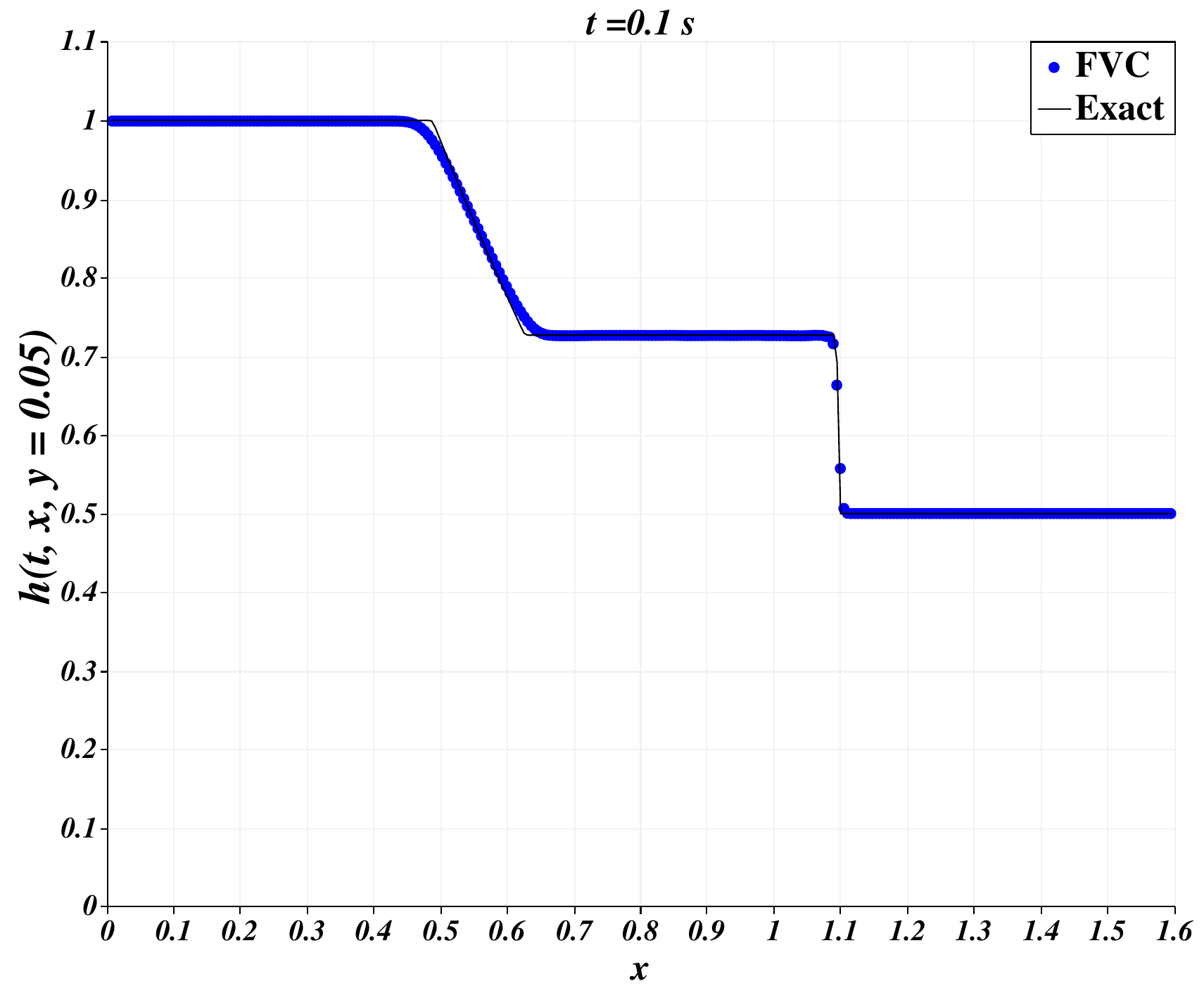} 
	\end{subfigure}
	\begin{subfigure}[b]{0.5\linewidth}
		\centering
		\includegraphics[width=.99\textwidth,height=.2\textheight]{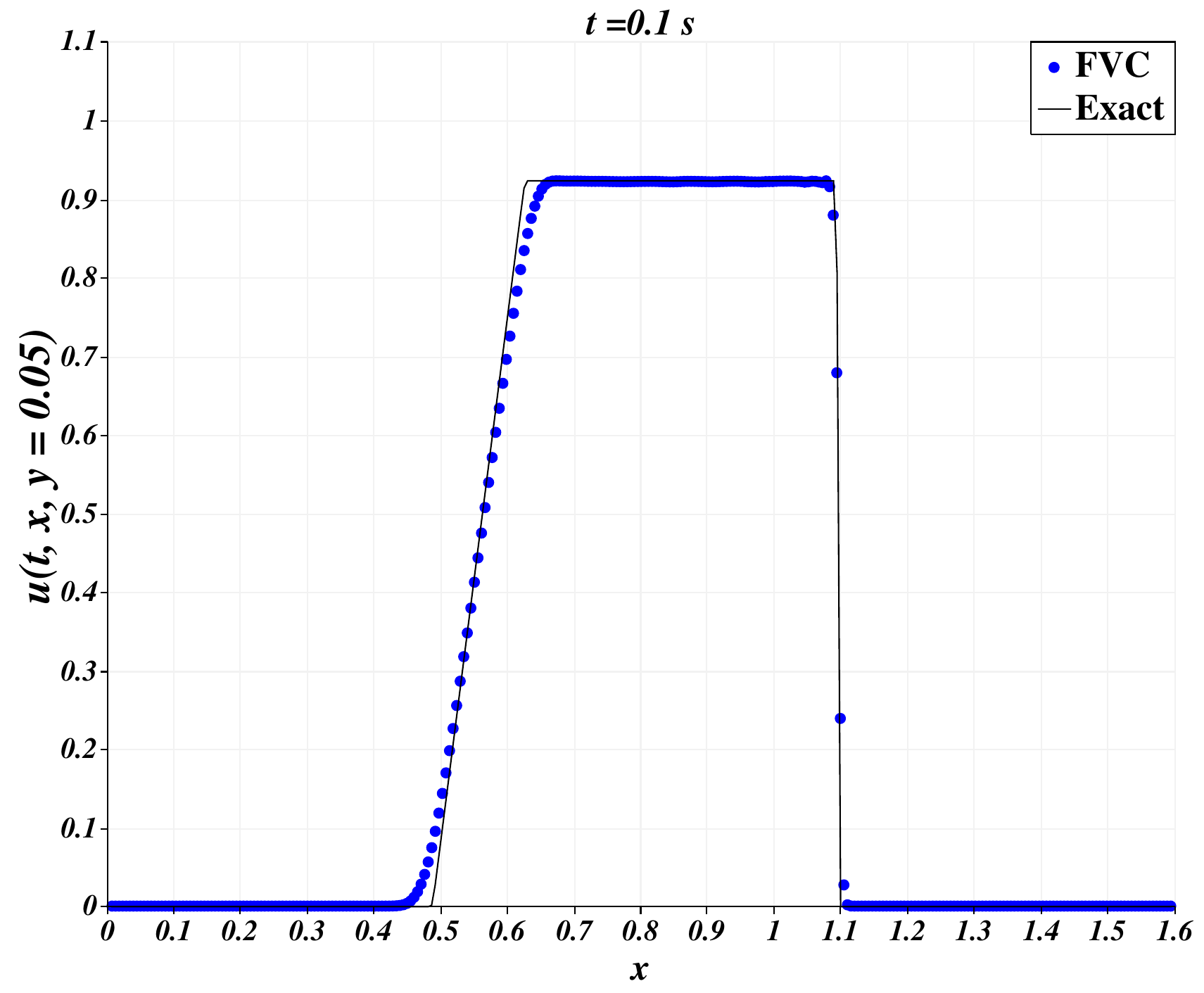} 
	\end{subfigure}
	\caption{Comparison of results for $h_r/h_l $ = 0.5 at t = 0.1s. 
		 Left: water height $h$.  right: longitudinal velocity $u$.}
\label{fig:dam05}
\end{figure}

The first simulation concerns a river flow with $h_r/h_l$ = 0.5. The \figurename{\ref{fig:dam05}} shows the cross-section at $y = 0.05$ of the evolution of the water depth and the longitudinal velocity. Excellent agreement is obtained between the numerical and analytical results, this is clear from Table \ref{Tab2} where the $L^1$ error and the accuracy order of $h$ and $hu$ are presented respectively. The comparison shows that, under this condition, our scheme, can accurately predict the shock wave without creating oscillations.

\begin{figure}[H] 
	\begin{subfigure}[b]{0.5\linewidth}
		\centering
		\includegraphics[width=.99\textwidth,height=.2\textheight]{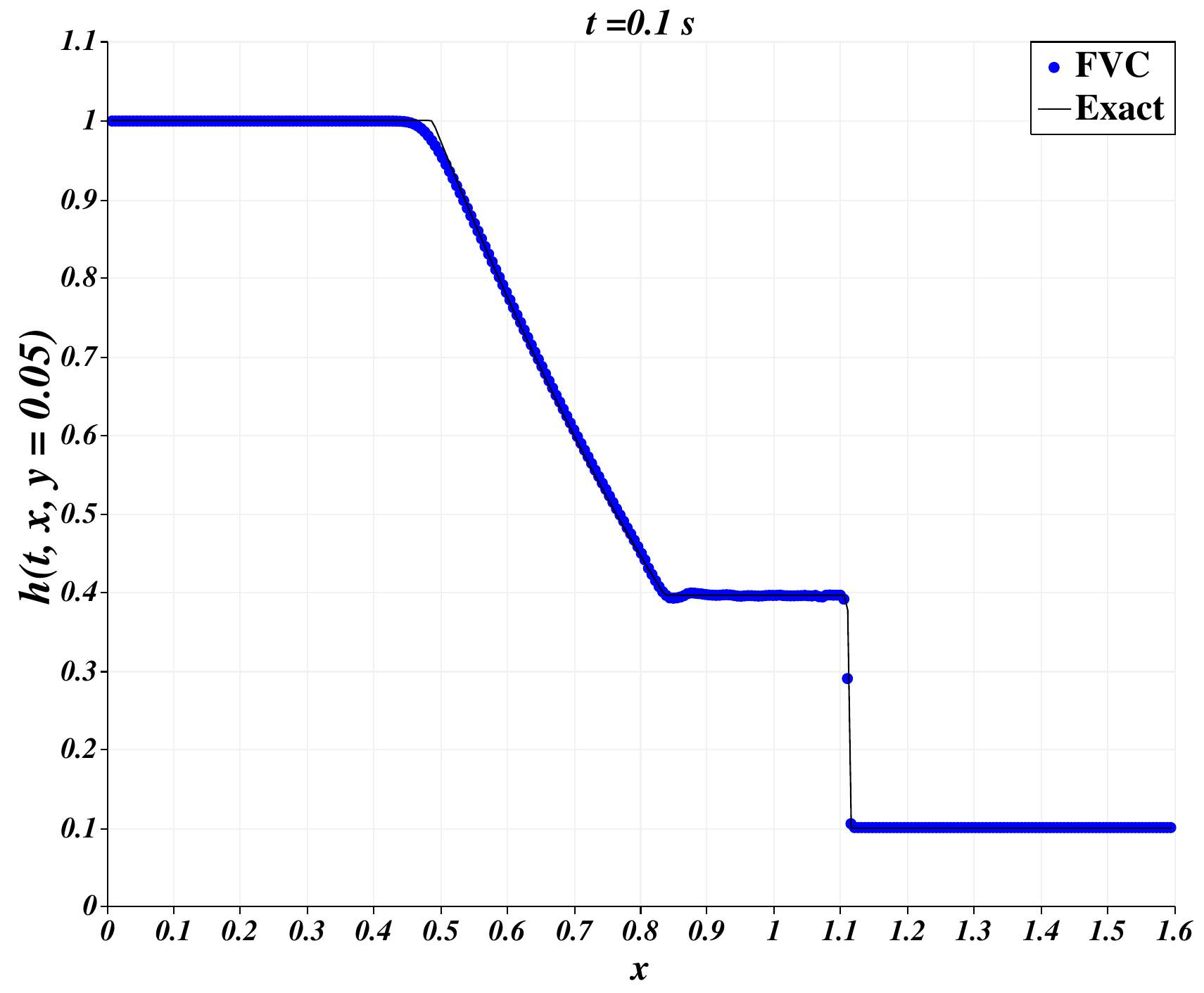} 
	\end{subfigure}
	\begin{subfigure}[b]{0.5\linewidth}
		\centering
		\includegraphics[width=.99\textwidth,height=.2\textheight]{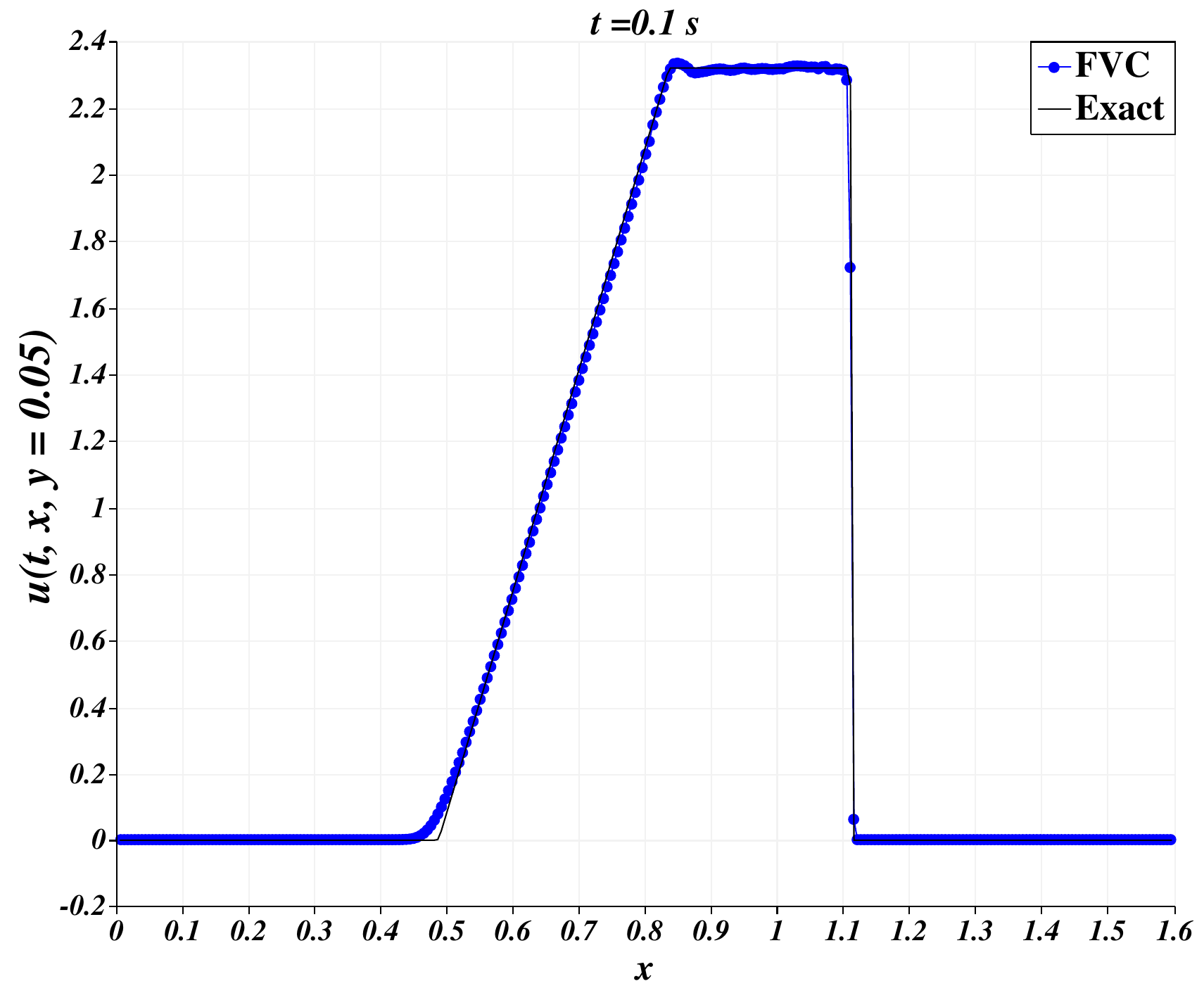} 
	\end{subfigure}
	\caption{Comparison of results for $h_r/h_l $ = 0.1 at t = 0.1s. 
		Left: water height $h$.  right: longitudinal velocity $u$.}
	\label{fig:dam01}
\end{figure}

\begin{figure}[H] 
	\begin{subfigure}[b]{0.5\linewidth}
		\centering
		\includegraphics[width=.99\textwidth,height=.2\textheight]{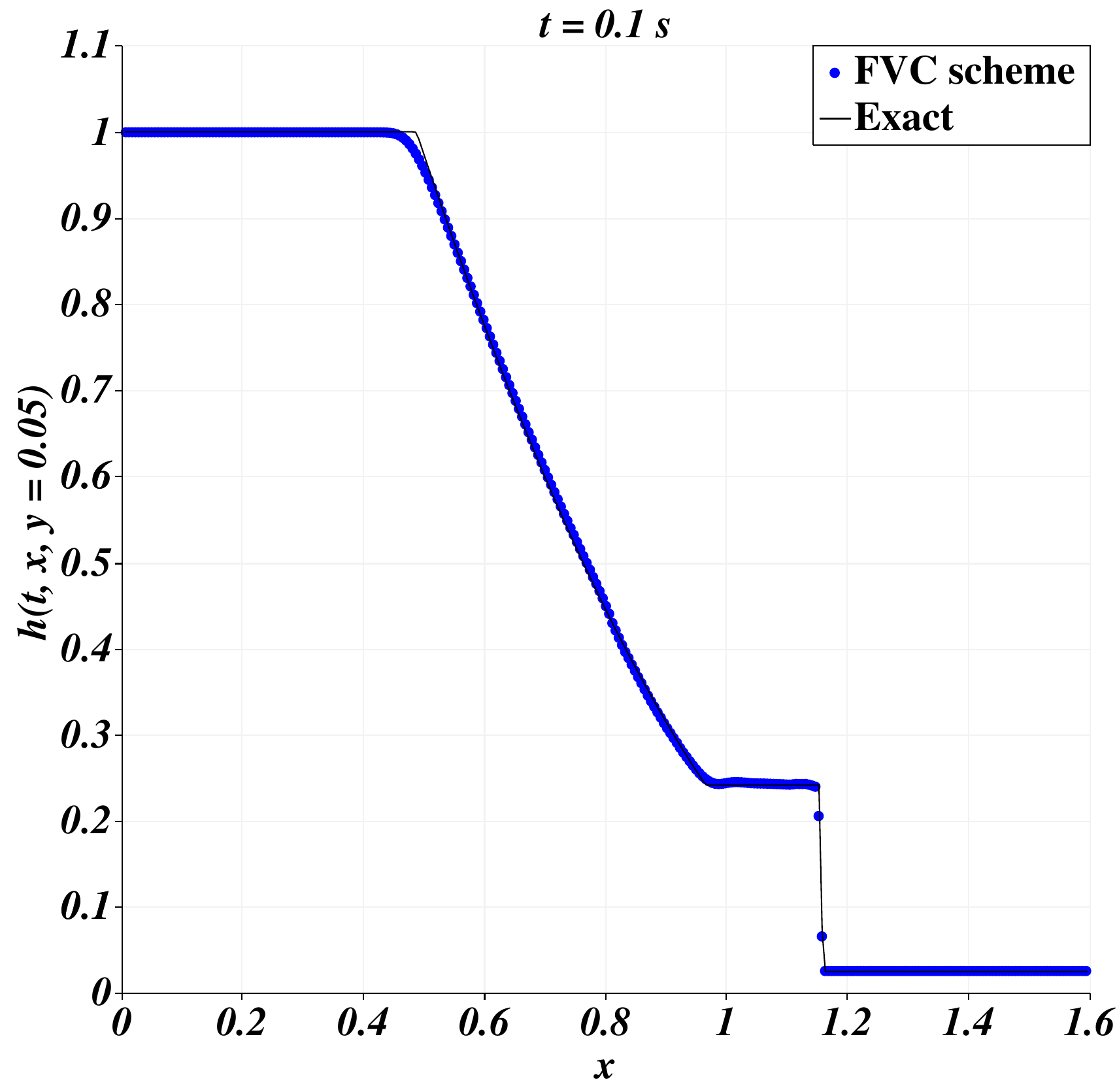} 
	\end{subfigure}
	\begin{subfigure}[b]{0.5\linewidth}
		\centering
		\includegraphics[width=.99\textwidth,height=.2\textheight]{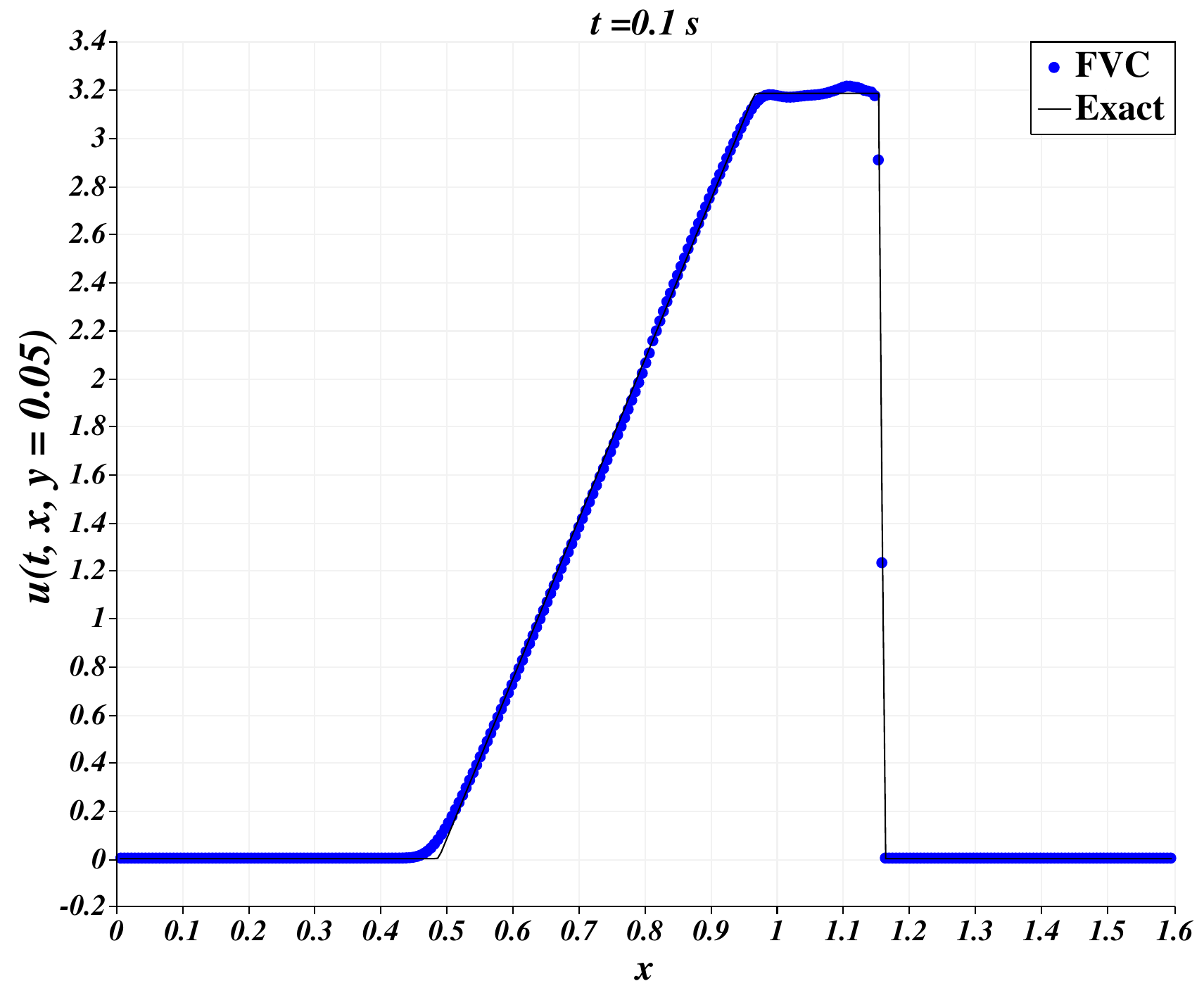} 
	\end{subfigure}
	\caption{Comparison of results for $h_r/h_l $ = 0.025 at t = 0.1s. 
		Left: water height $h$.  right: longitudinal velocity $u$.}
	\label{fig:dam025}
\end{figure}

In \cite{fennema1989implicit} the authors have proven that for a ratio $h_r/h_l$ smaller than $0.05$, most of the existing numerical models cannot give accurate results especially on the front. The last simulations, with $h_r/h_l$  = $0.1$ and $0.025$ (see \figurename{ \ref{fig:dam01}} and \figurename{\ref{fig:dam025}}), create supercritical flows downstream and sub-critical flows upstream.  When $h_r/h_l \leqslant 0.025$, a slight, non-physical oscillation occurs at the shock front.
\begin{table}[H]
	\caption{Relative $L^1$ error and CPU times  for dam break test at t = 0.1 s using FVC scheme on a different meshes.}\label{Tab2} 
	\centering
\begin{tabular}{cccccc}
	\hline 
  $\#$ Cells &\hspace{1.cm} Maximum of edges size &  \hspace{0.cm} Error in $h$  &  \hspace{1.cm} Error in $hu$    & \hspace{1.cm}  Ordre & \hspace{1.cm} CPU time (s)  \\ 
	\hline 
	5252& 0.0127  &2.612E-03  & \hspace{1.cm} 2.189E-02 & \hspace{1.cm}-  & 7.91\\ 
   10632& 0.00913 &1.650E-03  & \hspace{1.cm} 1.378E-02  & \hspace{1.cm}1.402  &10.66\\
   21224& 0.00666 & 1.045E-03  &\hspace{1.cm}  8.711E-03  & \hspace{1.cm}1.453 &19.96\\
   41776& 0.00479 & 6.255E-04  & \hspace{1.cm} 5.169E-03  & \hspace{1.cm}1.583 &30.21\\
	\hline 
\end{tabular} 
\end{table}
\vspace{-0.5cm}
The ratio $h_r/h_l $ is largely responsible for the problem of numerical instabilities that occur in the simulation of torrential flow due to dam failure. The difficulty of the problem increases with the decrease of the ratio $h_r/h_l$. We conclude that the proposed approach is well able to simulate torrential flows with a good capture of the shock front for the ratio $h_r/h_l $ very small.
\subsection{Tidal wave flow over an irregular bed}\label{irrebed}
Here we propose to study a tidal flow on a very irregular  bed. This test case has been proposed in several works to validate and to test the C-property of their approximation method \cite{zhou2001surface,benkhaldoun2010new}. So this test example is a good illustration of the significance of the source term treatment for practical applications to natural watercourses. It is well known that in the shallow water equations describing a flow over a very irregular background, the source terms become dominant and may cause undesirable numerical instabilities. Therefore this test case  allows us to test the reliability and robustness of the proposed model when the bottom variation is fast and irregular. The bed topography is defined in Appendix (Table of bed elevation $ Z(x)$ and its illustration ). The initial and boundary conditions are constructed from the asymptotic analytical solution is given by

$\displaystyle h(t,x) = h_0 + 4 - Z(x) - 4\sin(\pi (\frac{4t}{86400} + \frac{1}{2}))   $, \  \ \ \  $\displaystyle u(t,x) = \frac{\pi(x  - 1500)}{5400\cdot h(t,x)} \cos(\pi(\frac{4t}{86400} + \frac{1}{2})) $ \\

but with,  \  \ \ \ $\displaystyle h_0 = 16 \ m  $, \ \ \ $\displaystyle h(0,x) = h_0 - Z(x). $ 

\begin{figure}[H] 
	\begin{subfigure}[b]{0.5\linewidth}
		\centering
		\includegraphics[width=.99\textwidth,height=.22\textheight]{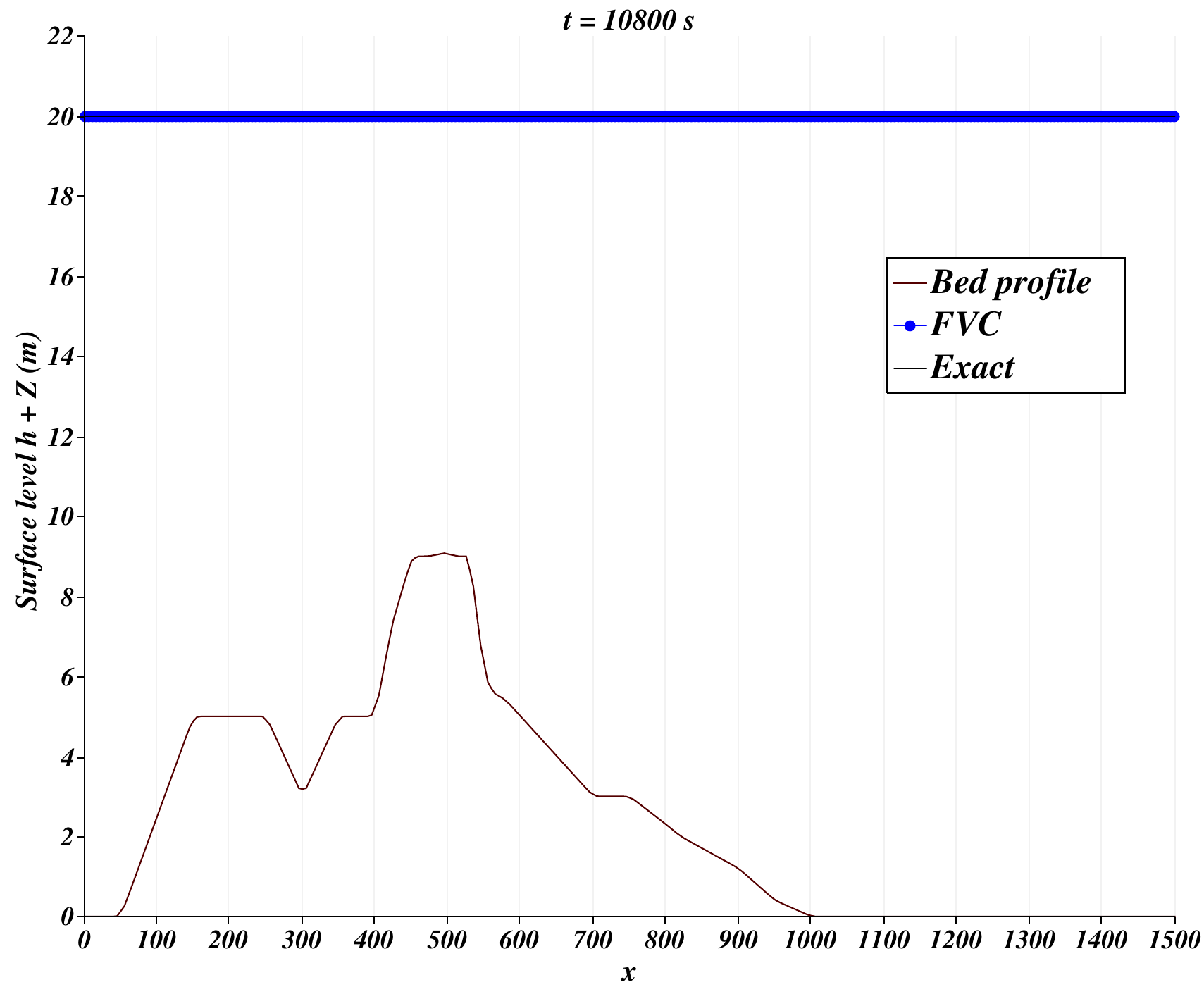} 
	\end{subfigure}
	\begin{subfigure}[b]{0.5\linewidth}
		\centering
		\includegraphics[width=.99\textwidth,height=.22\textheight]{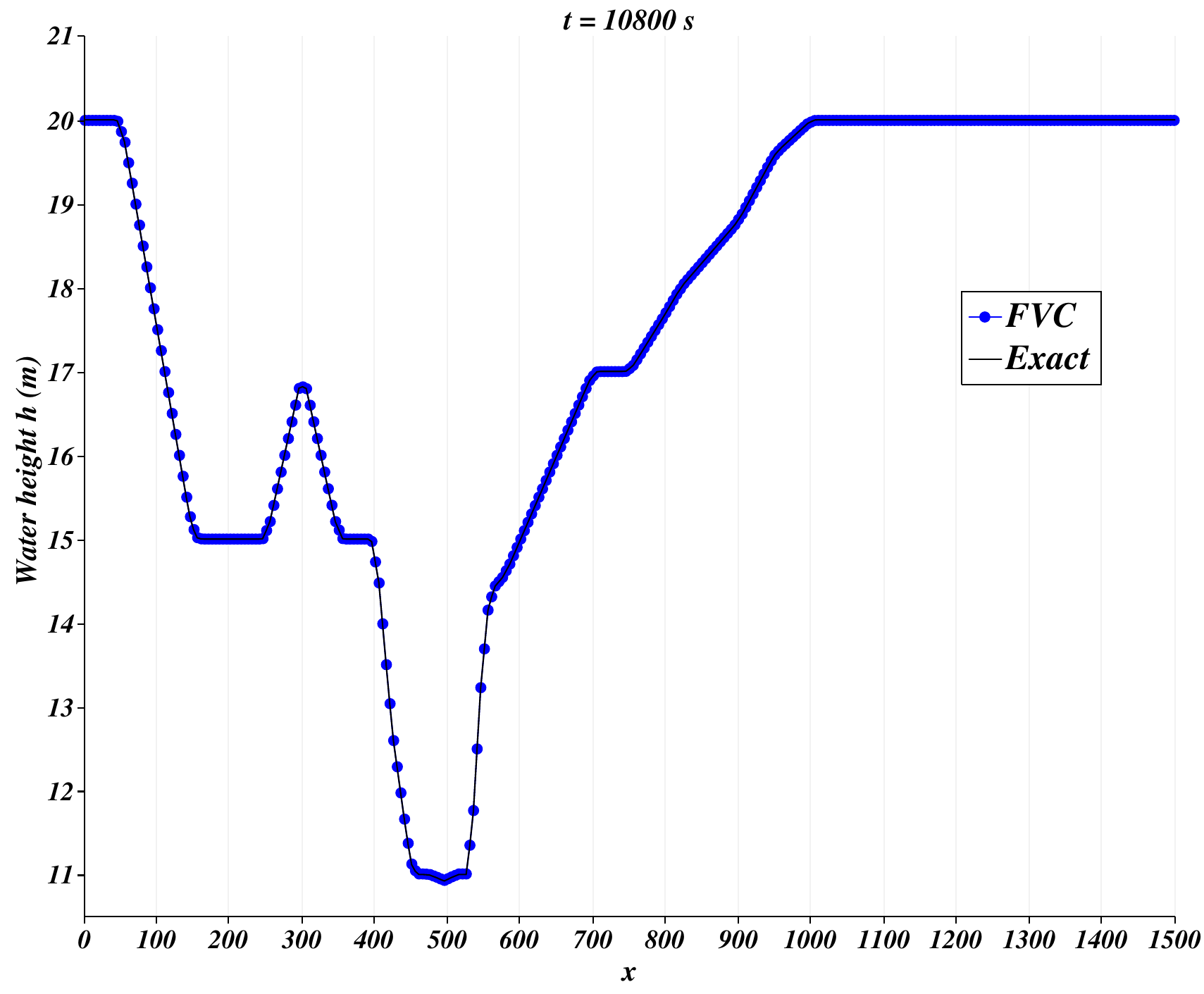} 
	\end{subfigure}
	\caption{ Tidal wave flow over an irregular bed. Left: Comparison of surfaces level $h + Z$ at $t = 10800 \ s$ and bed $ Z $. Right: Comparison of water height $h$.  }
	\label{fig:tidalhpz}
\end{figure}

\begin{figure}[H] 
	\begin{subfigure}[b]{0.5\linewidth}
		\centering
		\includegraphics[width=.99\textwidth,height=.22\textheight]{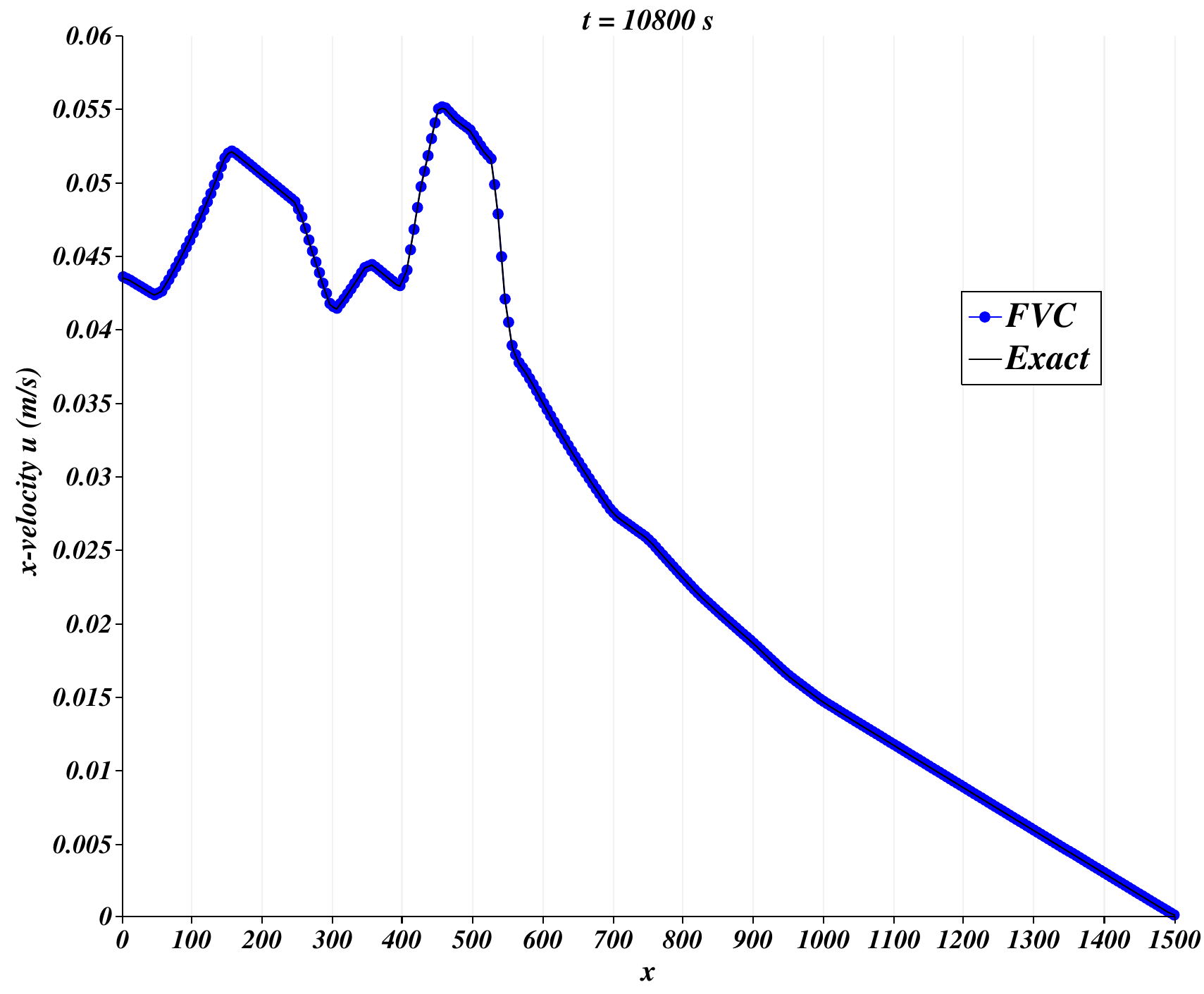} 
	\end{subfigure}
	\begin{subfigure}[b]{0.5\linewidth}
		\centering
		\includegraphics[width=.99\textwidth,height=.22\textheight]{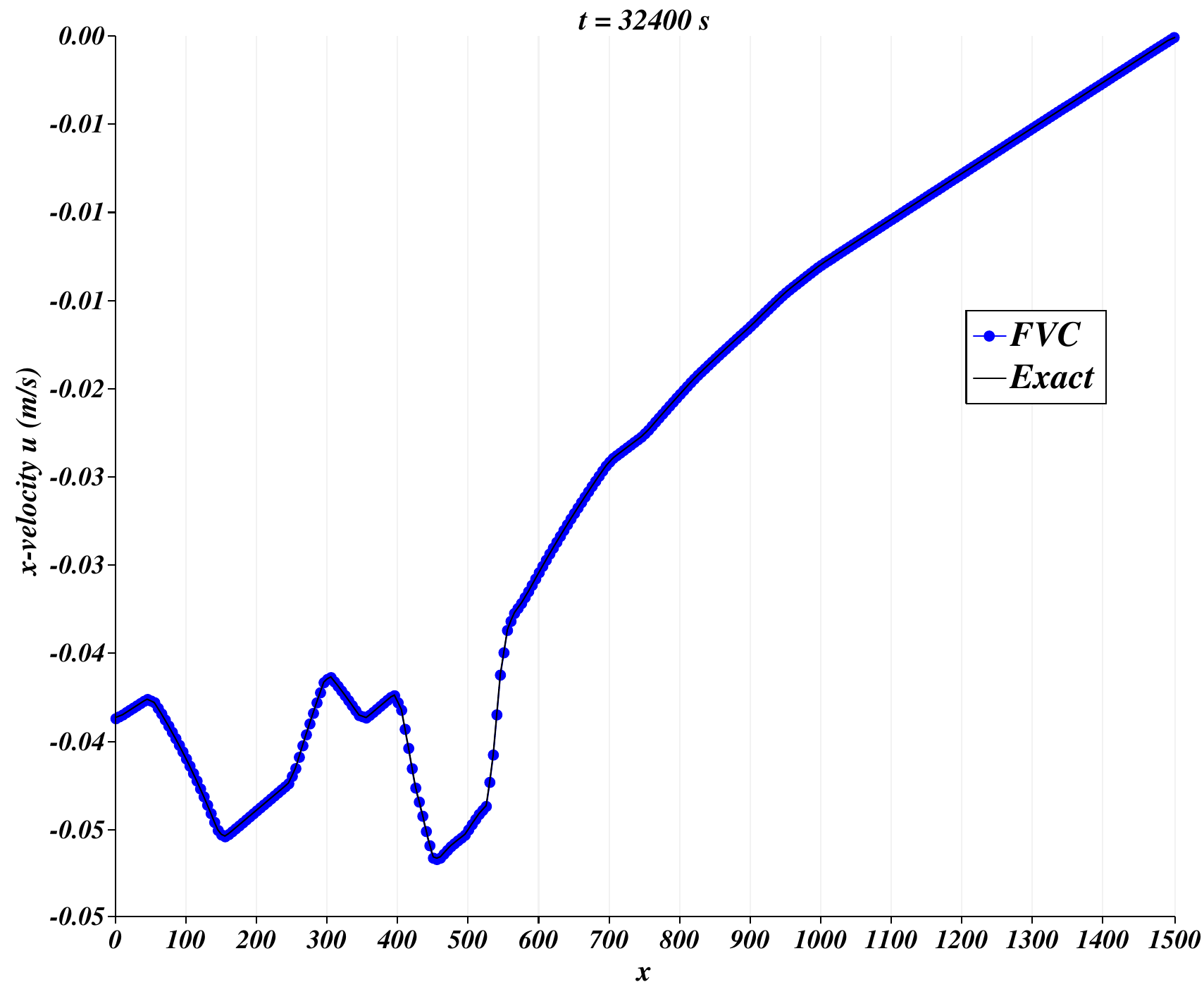} 
	\end{subfigure}
	\caption{ Comparison of water velocity $u$ in tidal wave flow over an irregular bed.  Left: at $t = 10800 \ s$.  Right: at $ 32400 \ s $.  }
	\label{fig:tidalu}
\end{figure}

In order to compare the numerical results of our approach with the analytical solution, we choose two results, at $t = 10800 \ s$ and at $ t = 32400 \ s$.
In \figurename{\ref{fig:tidalhpz}} we present a comparison between the approach surface level and the analytical solution at $t = 10800\  s $ as well as the water height at the same time using a mesh of 200 grid points in x-direction. We also include in \figurename{\ref{fig:tidalu}} a comparison between the water velocity generated by the FVC scheme and the analytical velocity at $t = 10800 \ s$  then at  $t = 32400 \ s$.
In Table \ref{Tab3}  we present a comparison between the exact solution and the solution generated by FVC scheme using the relative L$^1$ error. 
An excellent agreement is obtained between the numerical and analytical solutions. This confirms that the proposed scheme is also accurate for tidal flow over an irregular bed.  Moreover, these results are qualitatively in good agreement with those published in \cite{zhou2001surface,benkhaldoun2010new}.
\begin{minipage}[H]{0.4\linewidth}
 
The accuracy of the proposed scheme in the treatment of the source terms has been identified. The numerical errors produced by the model remain very low despite the fact that the mesh is coarse and the bed is very irregular. The results of this test case confirm the good performance of our approach in the treatment of the source terms, while avoiding undesirable numerical errors due to the rapid variation of the bed.
\end{minipage}\hspace{.5cm}
\begin{minipage}[H]{0.55\linewidth}
	\vspace{-1.cm}
\begin{table}[H]
	\caption{Relative $L^1$ error and CPU times  for the tidal wave flow over an irregular bed using FVC scheme.}\label{Tab3} 
	\centering
	\begin{tabular}{ccccc}
		\hline\\[-0.5em] 
		$t_{end}$ &Error in $h$ &  Error in $hu$  & Error in $h + Z$  & CPU time (s)  \\ 
		\hline\\[-0.5em] 
	    10800&  1.254E-05  &4.230E-03   & 1.012E-05   & 545.95\\ 
		32400&  7.719E-06 & 5.731E-03  & 1.241E-05   &1583.41\\[0.5em] 
		\hline 
	\end{tabular} 
\end{table}
\end{minipage}

\subsection{ Flow over a non-flat irregular  bed }\label{Cprop}


We consider the example of water flow in a two-dimensional
channel including an irregular bed, a similar test has been proposed in \cite{benkhaldoun2010new}. The mathematical formulation consists
of solving the shallow water system (\ref{1}) without Coriolis force and subjected to Neumann boundary conditions.
The initial conditions as follow \\
$\displaystyle h(0,x,y) = 1 - Z(x,y) \ m, \quad \quad  \displaystyle u(0,x,y) = v(0,x,y) = 0\  m/s $, \\
where the bed profile is defined by: \ \ 
$\displaystyle Z(x,y)=\sum_{k=1}^{5} a_k\exp(-\frac{(x-x_k)^2 + (y-y_k)^2}{\sigma_k^2})$,\\
with \ \  $(a_1, \sigma_1^2, x_1, y_1)=(0.75, 2, -4,5)$,\  \ $(a_2, \sigma_2^2,x_2,y_2 ) = (0.7,2, -2.5,2.5) $,\  \  $ (a_3, \sigma_3^2,x_3,y_3 ) = (0.65,3.3,0,0) $,\\[.2cm]
 $ (a_4, \sigma_4^2, x_4,y_4) = (0.6,2.5,3,-2) $, and $ (a_5, \sigma_5^2,x_5,y_5) = (0.55,1.48,5,-4) $.\\
The purpose of this test example is to verify the achievement of the
C-property for the FVC scheme applied to shallow water flows over non-flat bed.
\begin{figure}[H] 
	\begin{subfigure}[b]{0.5\linewidth}
		\centering
		\includegraphics[width=.99\textwidth,height=.22\textheight]{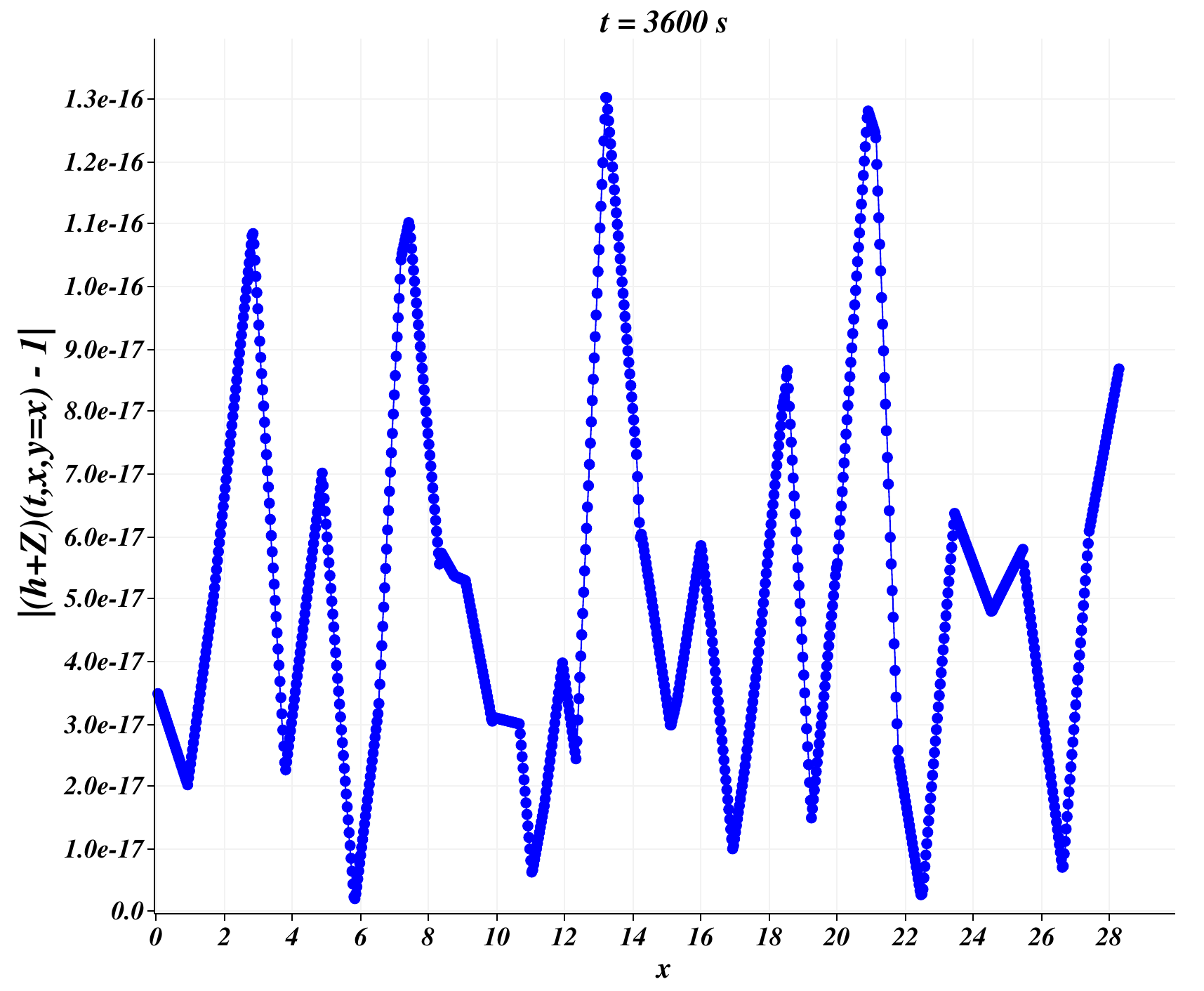}
	\end{subfigure}
	\begin{subfigure}[b]{0.5\linewidth}
		\centering
		\includegraphics[width=.99\textwidth,height=.22\textheight]{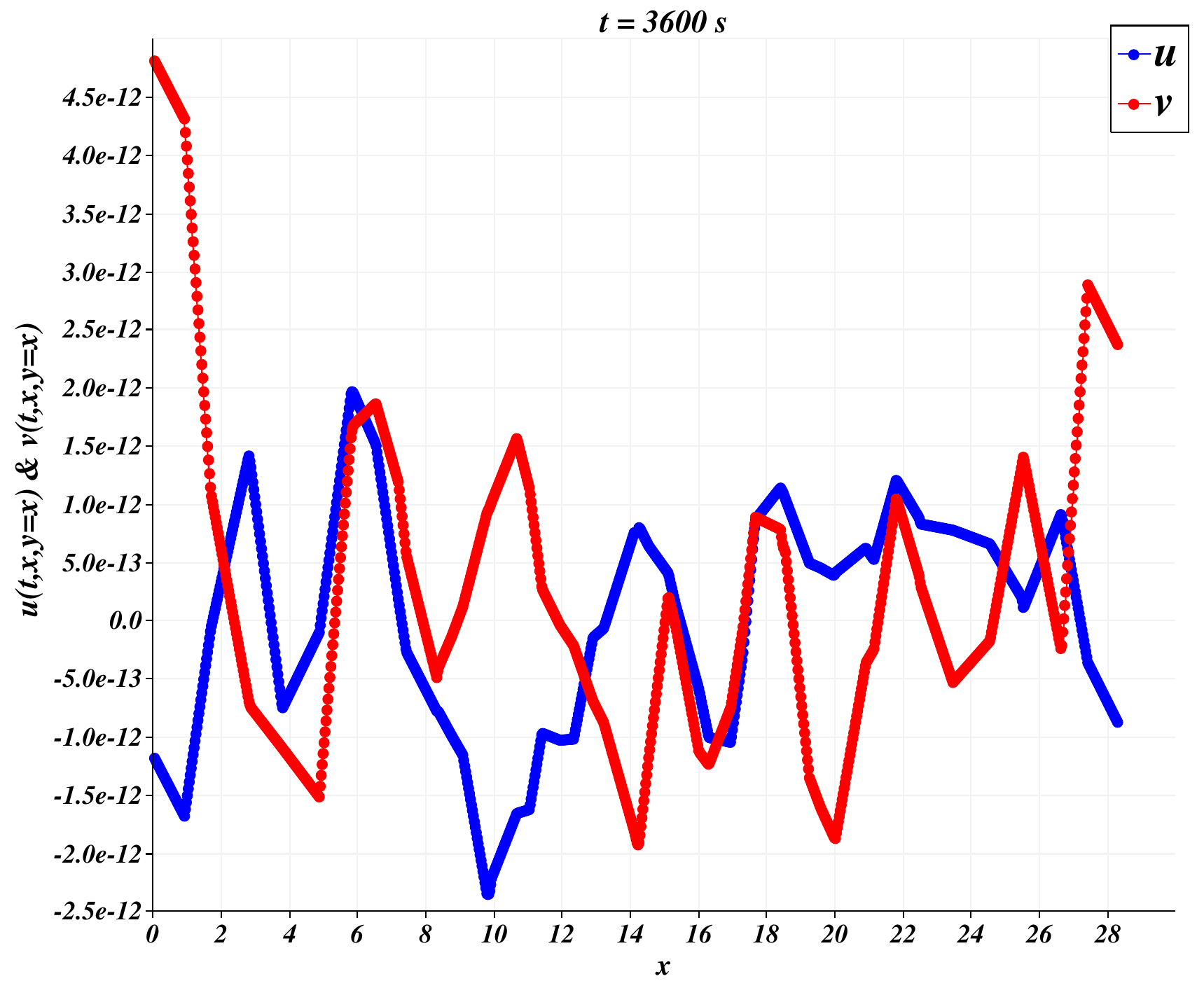} 
	\end{subfigure}
	\caption{Cross-sections at $y = x$. Left: the absolute error of the free-surface for the lake at rest. Right: velocities values after 1 hour.}
	\label{fig:hpz_uv}
\end{figure}
The C-property of well-balanced FVC scheme on unstructured meshes is checked in this example. 
As expected the water free-surface remains constant during the simulation time. The velocities and the error presented in \figurename{\ref{fig:hpz_uv}} also show that the equilibrium of the lake at rest is verified. All this shows that the proposed  FVC scheme perfectly preserves the C-property.
\subsection{ Circular dam-break problem}\label{Circular}
We consider the benchmark problem proposed in \cite{kuo2000nonlinear} to study cyclone/anticyclone asymmetry in nonlinear geostrophic adjustment. We solve the shallow water system (\ref{1}) with a Coriolis effect on a
non-flat bottom in the spatial domain $ \Omega = [-10, 10] \times [-10, 10] $ subjected
to Neumann boundary conditions and the following initial conditions\\

$\displaystyle h(0,x,y) = 1 +\frac{1}{4}\left(1-\tanh \left( \frac{\sqrt{ax^2+by^2}-1}{c}  \right)  \right)$,\quad
$ \displaystyle u(0,x,y) = v(0,x,y) = 0 m/s,$\\[0.5em]
where $\  a = \frac{5}{2},$ \ \ \ $b =\frac{2}{5},$ \ \ $c = 0.1,$  \ \ \ $f_c = 1\  Kg.m/s ^2$ \ and \ \ \ $ g = 1\  m/s ^2$.\\

The bottom profile has the following expression: \ $\displaystyle Z(x,y)= 0.3\left(1 + \tanh(\frac{3x}{2})\right)$.\\

Let's start by looking at the behaviour of this phenomenon in a domain with a flat bottom.  
The \figurename{\ref{fig:Circular1} shows the representation of the water level calculated at different times for this test case with $Z(x,y) = 0$. As can be seen, a hole has formed and water is flowing out of the deepest region as a rarefaction wave progresses outwards. It is clear from the results presented that the initial elliptical mass imbalance evolves in a non-axisymmetric manner.  The two expected shock waves are very well captured by the proposed FVC method. These results are qualitatively in good agreement with those published in \cite{benkhaldoun2015projection,kuo2000nonlinear}. In \figurename{\ref{fig:Circular2}} we exhibit the results for the velocity field corresponding to the plots \figurename{\ref{fig:Circular1}}.  As can be seen the two shock waves originated behind the water elevation are slowly spinning clockwise in the computational domain. The velocity field is well represented by the FVC method and re-circulation regions within the flow domain are well captured.
\begin{figure}[H] 
	\begin{subfigure}[b]{0.333\linewidth}
		\centering
		\includegraphics[width=1.1\textwidth,height=.18\textheight]{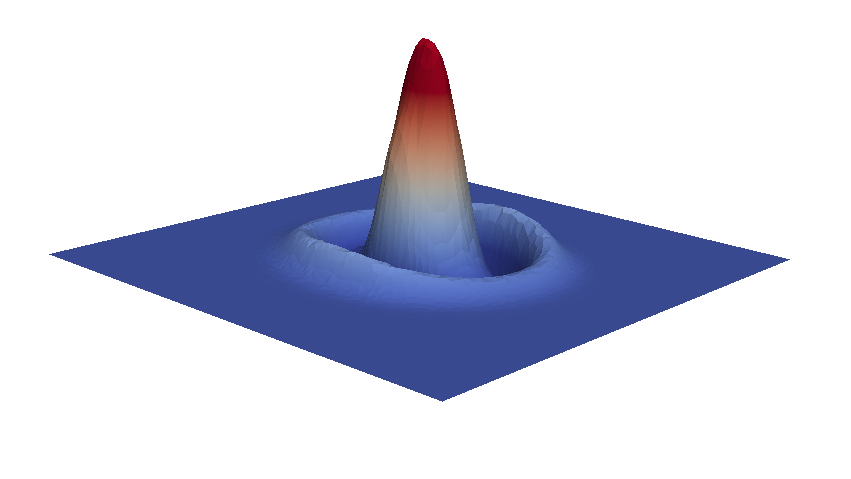} 
	\end{subfigure}
	\begin{subfigure}[b]{0.333\linewidth}
		\centering
		\includegraphics[width=1.1\textwidth,height=.18\textheight]{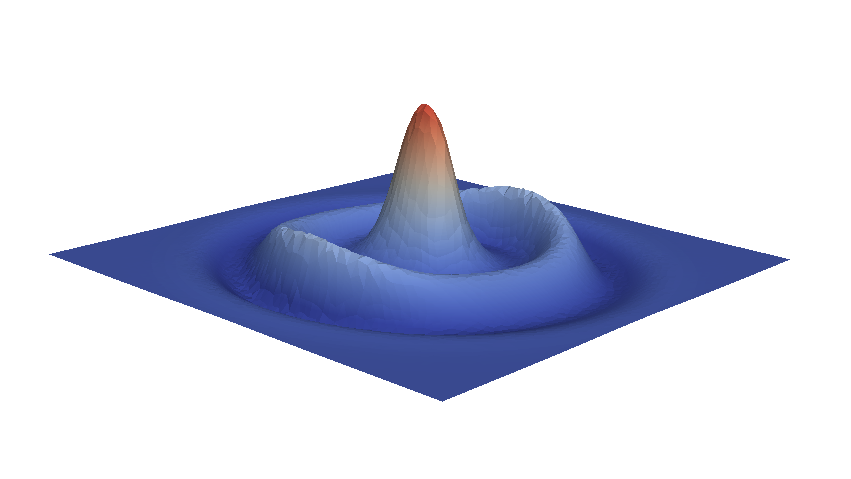} 
	\end{subfigure}
	\begin{subfigure}[b]{0.333\linewidth}
		\centering
		\includegraphics[width=1.1\textwidth,height=.18\textheight]{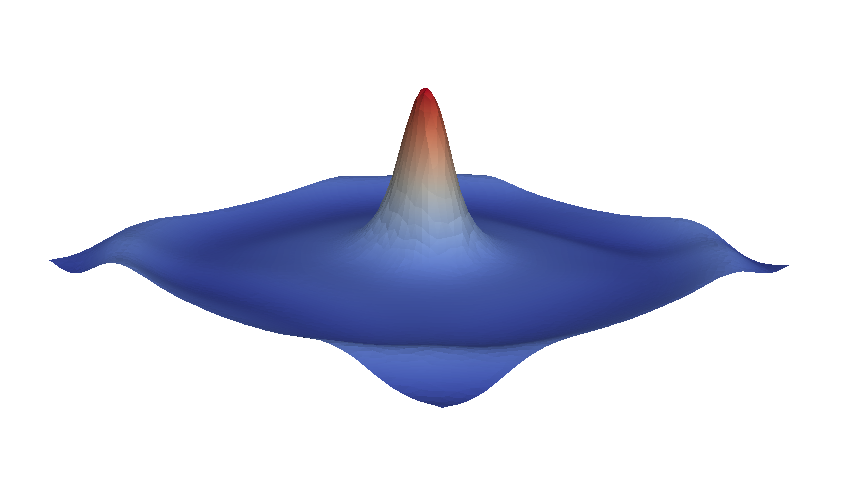} 
	\end{subfigure}
\caption{Water depth for the circular dam-break problem on flat bottom obtained at different times using a mesh with $10040$ cells. From top to bottom $t = 4$ s,  $8$ s and $16$ s.}
\label{fig:Circular1}
\end{figure}
\vspace{-.5cm}
\begin{figure}[H] 
	\begin{subfigure}[b]{0.33\linewidth}
		\centering
		\includegraphics[width=1.07\textwidth,height=.19\textheight]{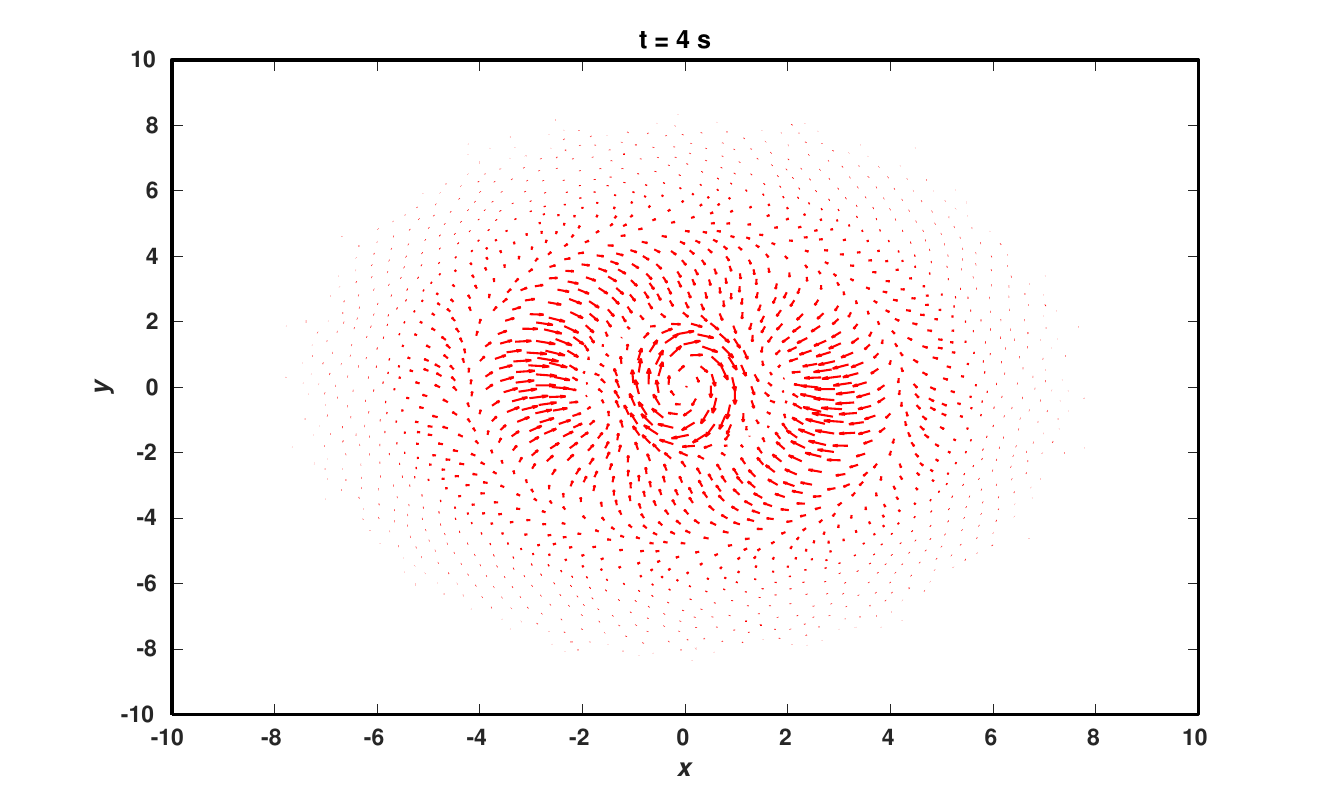} 
	\end{subfigure}
	\begin{subfigure}[b]{0.33\linewidth}
		\centering
		\includegraphics[width=1.07\textwidth,height=.19\textheight]{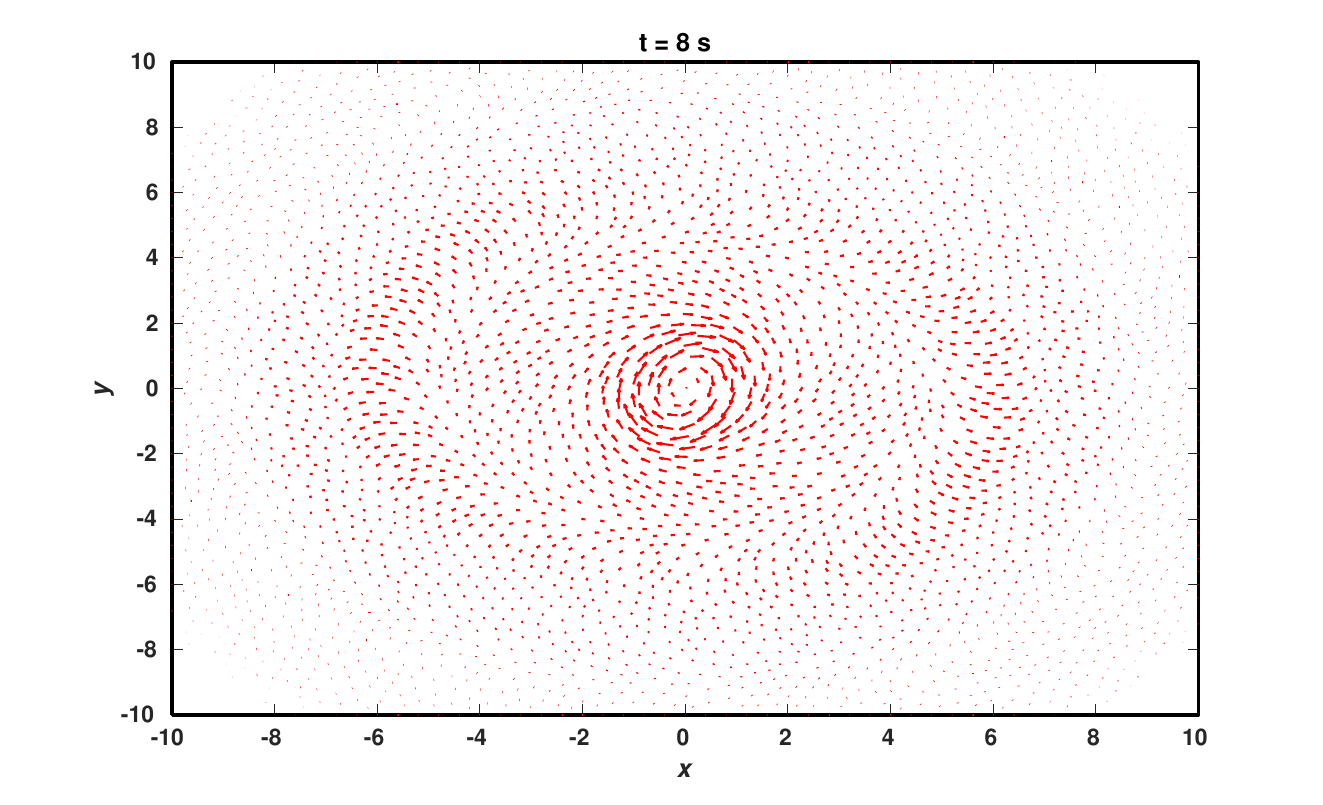} 
	\end{subfigure}
	\begin{subfigure}[b]{0.33\linewidth}
		\centering
		\includegraphics[width=1.07\textwidth,height=.19\textheight]{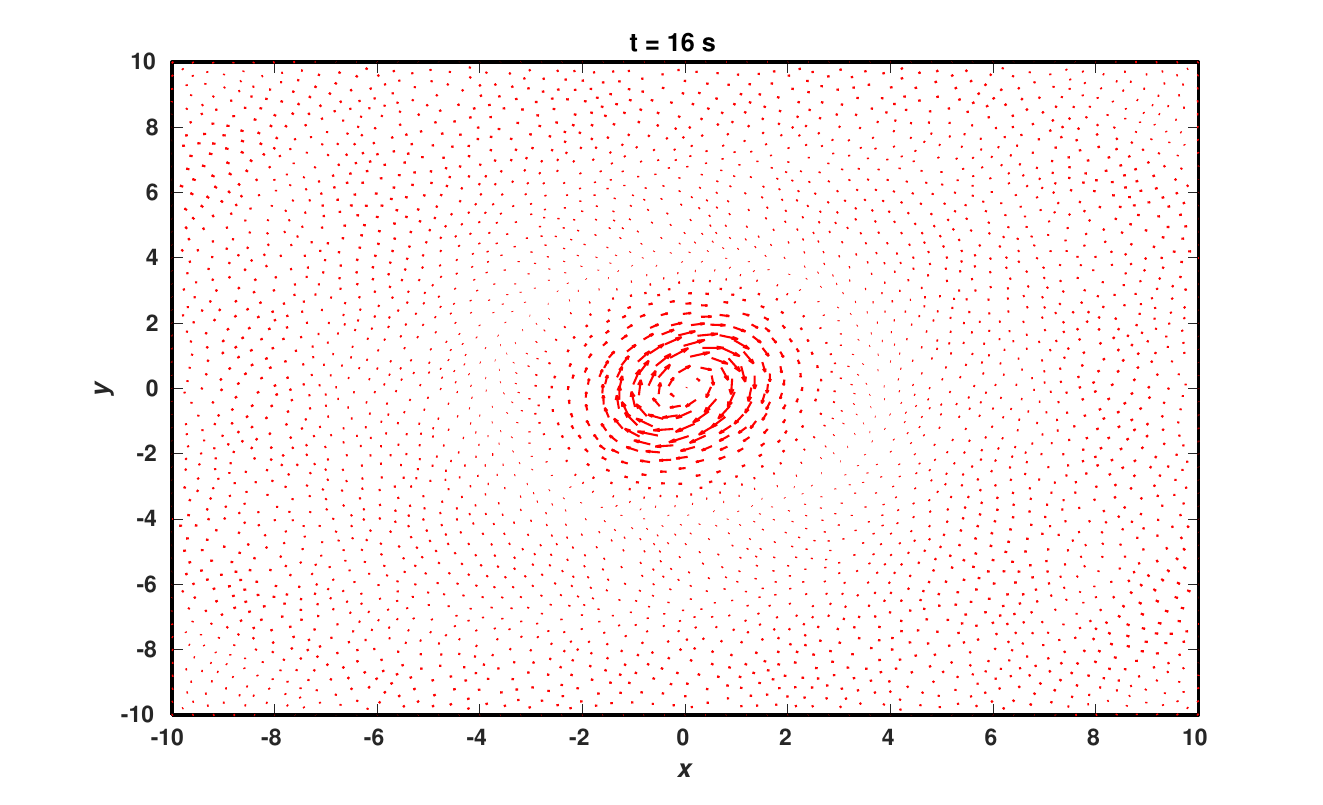} 
	\end{subfigure}
\caption{Velocity fields for the circular dam-break problem corresponding to the plots represented in \figurename{\ref{fig:Circular1}}.}
\label{fig:Circular2}
\end{figure}
Let's move on to the case of the non-flat bottom to assess the performance of our FVC scheme on unstructured meshes to solve the circular dam-break problem on a non-flat bottom. The \figurename{\ref{fig:Circular3}} shows the calculated results for the water depth at $t = 2$ s, $8 $ s and $16 $ s using two meshes of 10040 and 40146 cells. The corresponding results for  velocity field are presented in \figurename{\ref{fig:Circular4}}. From a numerical point of view this test example is more difficult than the previous one as the flow is expected to exhibit complex features due to the interaction between the water surface and the bed. As in the previous test a hole has formed and the water drains from the deepest region as a rarefaction wave progresses outwards. However, a slower propagation is detected for the water free-surface in this test compared to the simulations on flat-bottom.
\begin{figure}[H] 
	\begin{subfigure}[b]{0.333\linewidth}
		\centering
		\includegraphics[width=1.1\textwidth,height=.19\textheight]{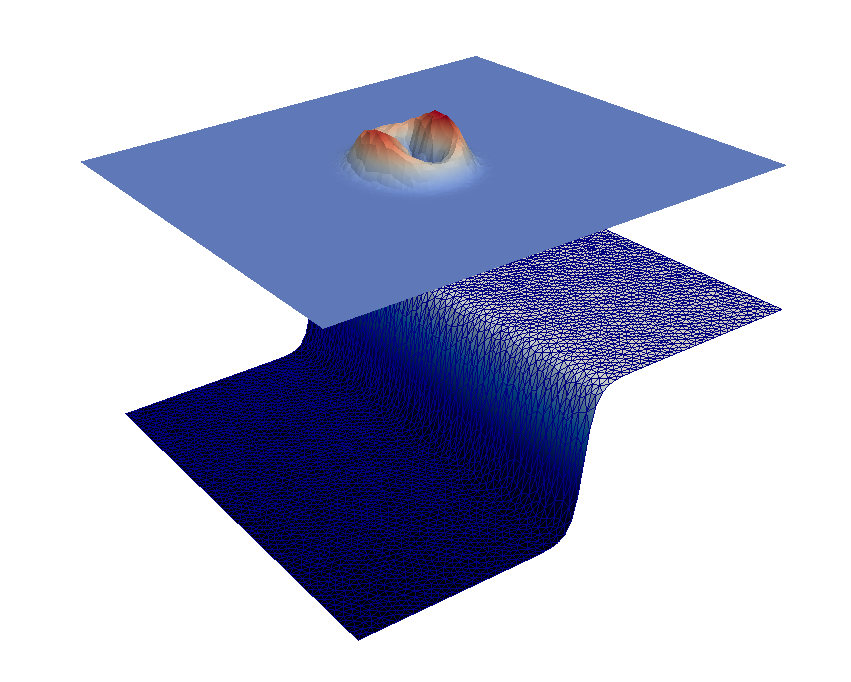} 
	\end{subfigure}
	\begin{subfigure}[b]{0.333\linewidth}
		\centering
		\includegraphics[width=1.1\textwidth,height=.19\textheight]{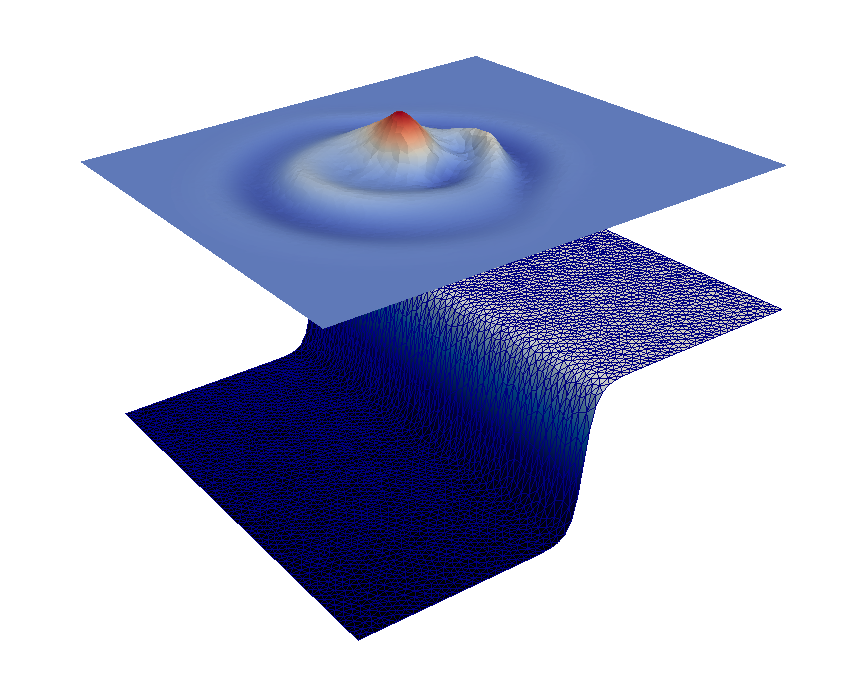} 
	\end{subfigure}
  	\begin{subfigure}[b]{0.333\linewidth}
		\centering
		\includegraphics[width=1.1\textwidth,height=.19\textheight]{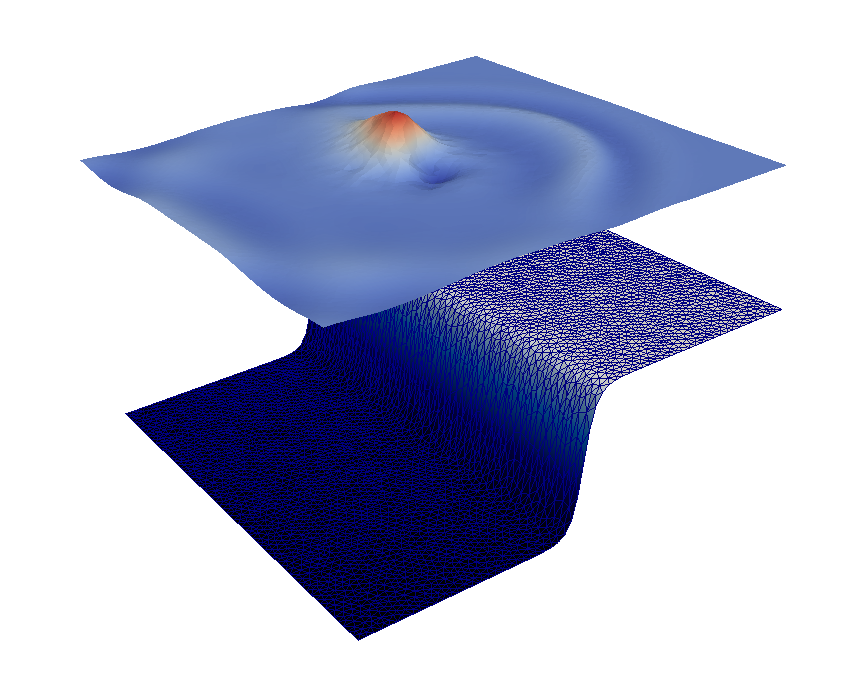} 
	\end{subfigure}
	\begin{subfigure}[b]{0.333\linewidth}
	\centering
	\includegraphics[width=1.1\textwidth,height=.19\textheight]{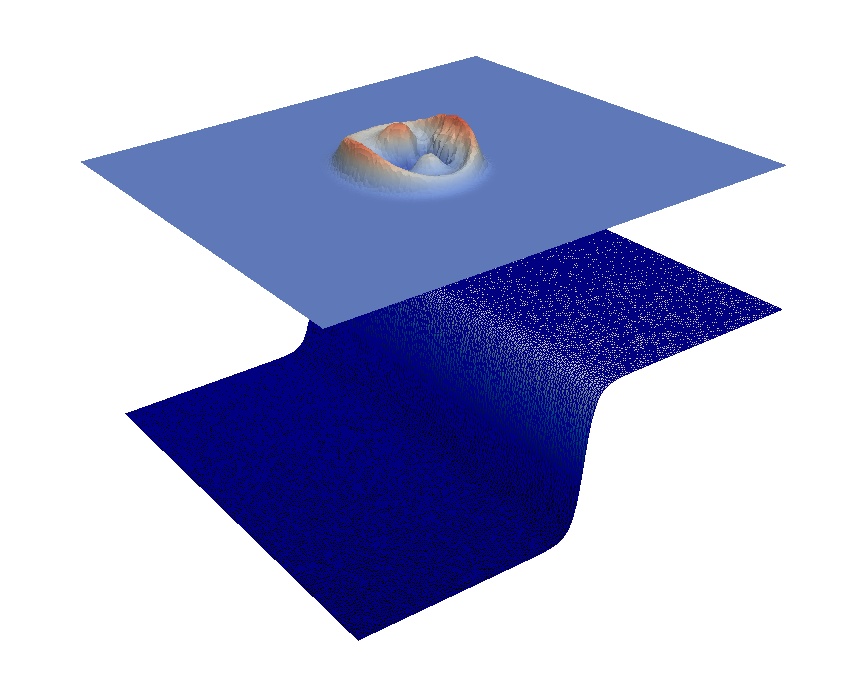} 
\end{subfigure}
\begin{subfigure}[b]{0.33\linewidth}
	\centering
	\includegraphics[width=1.1\textwidth,height=.19\textheight]{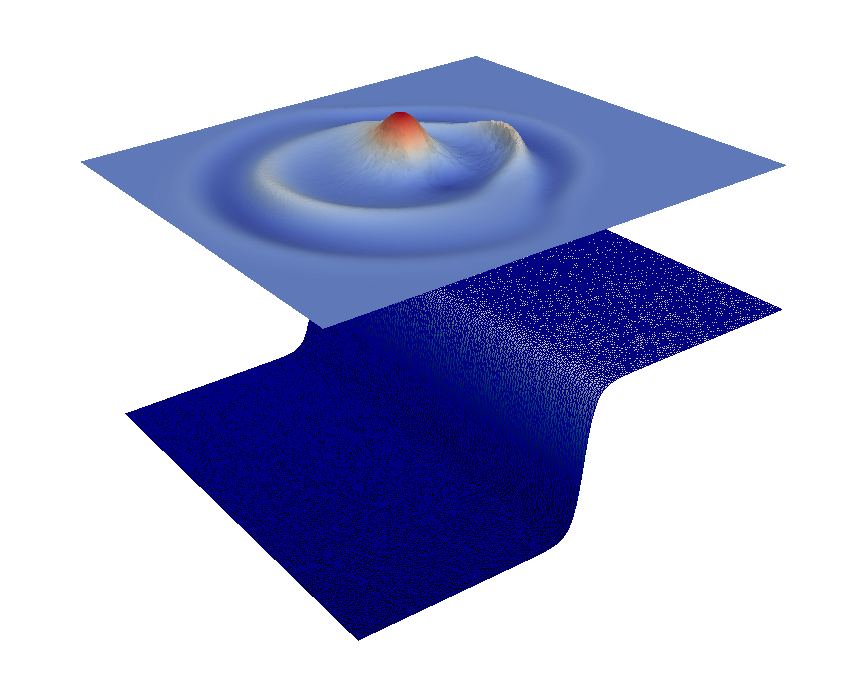} 
\end{subfigure}
\begin{subfigure}[b]{0.333\linewidth}
	\centering
	\includegraphics[width=1.1\textwidth,height=.19\textheight]{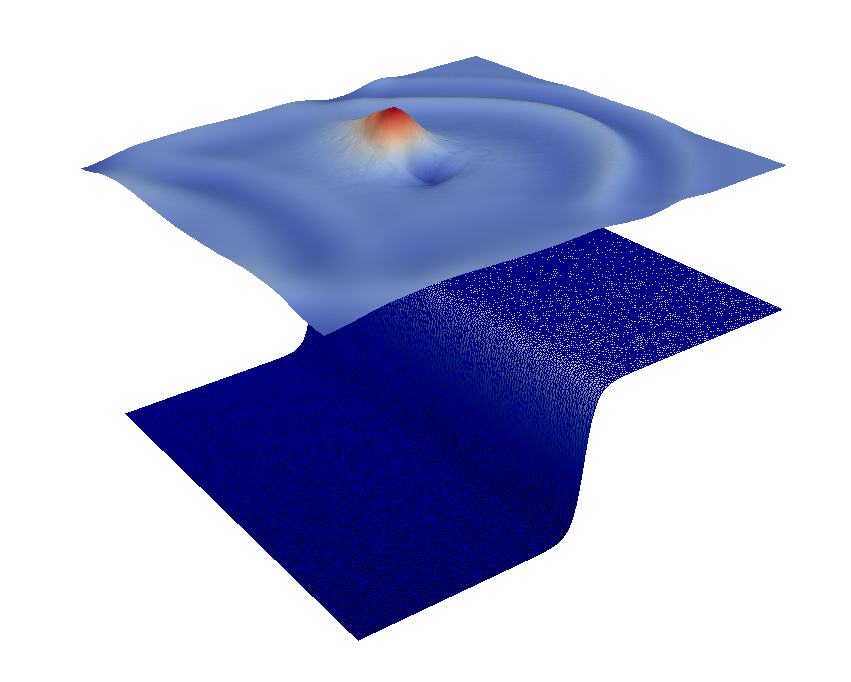} 
\end{subfigure}
	\caption{Water depth for the circular dam-break problem on non-flat bottom obtained using a mesh with 10040 cells (first row) and 40146 (second row). From left to right $t = 2$ s, $8$ s and $16$ s.}
	\label{fig:Circular3}
\end{figure}%
\vspace{-.4cm}
\begin{figure}[H] 
	\begin{subfigure}[b]{0.33\linewidth}
		\centering
		\includegraphics[width=1.05\textwidth,height=.19\textheight]{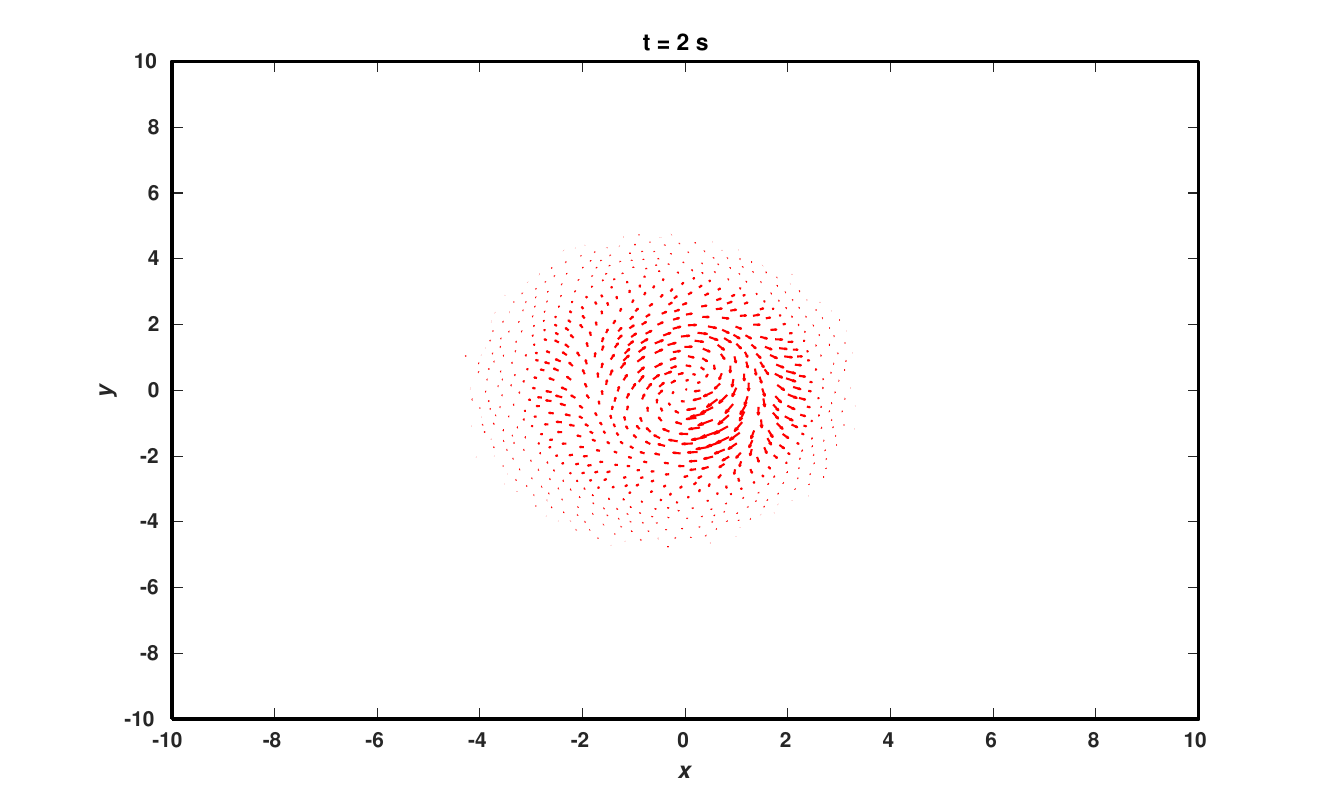} 
	\end{subfigure}
	\begin{subfigure}[b]{0.33\linewidth}
		\centering
		\includegraphics[width=1.05\textwidth,height=.19\textheight]{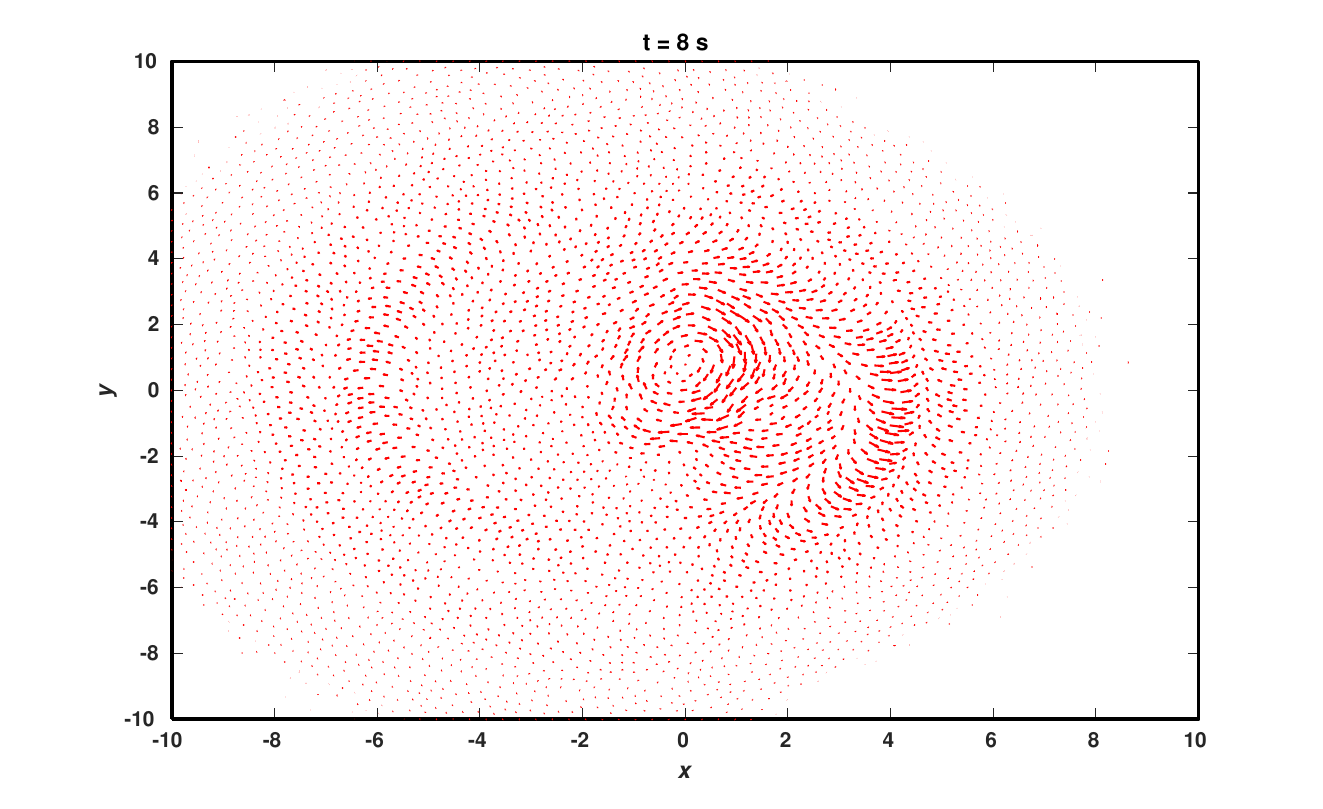} 
	\end{subfigure}
	\begin{subfigure}[b]{0.33\linewidth}
		\centering
		\includegraphics[width=1.05\textwidth,height=.19\textheight]{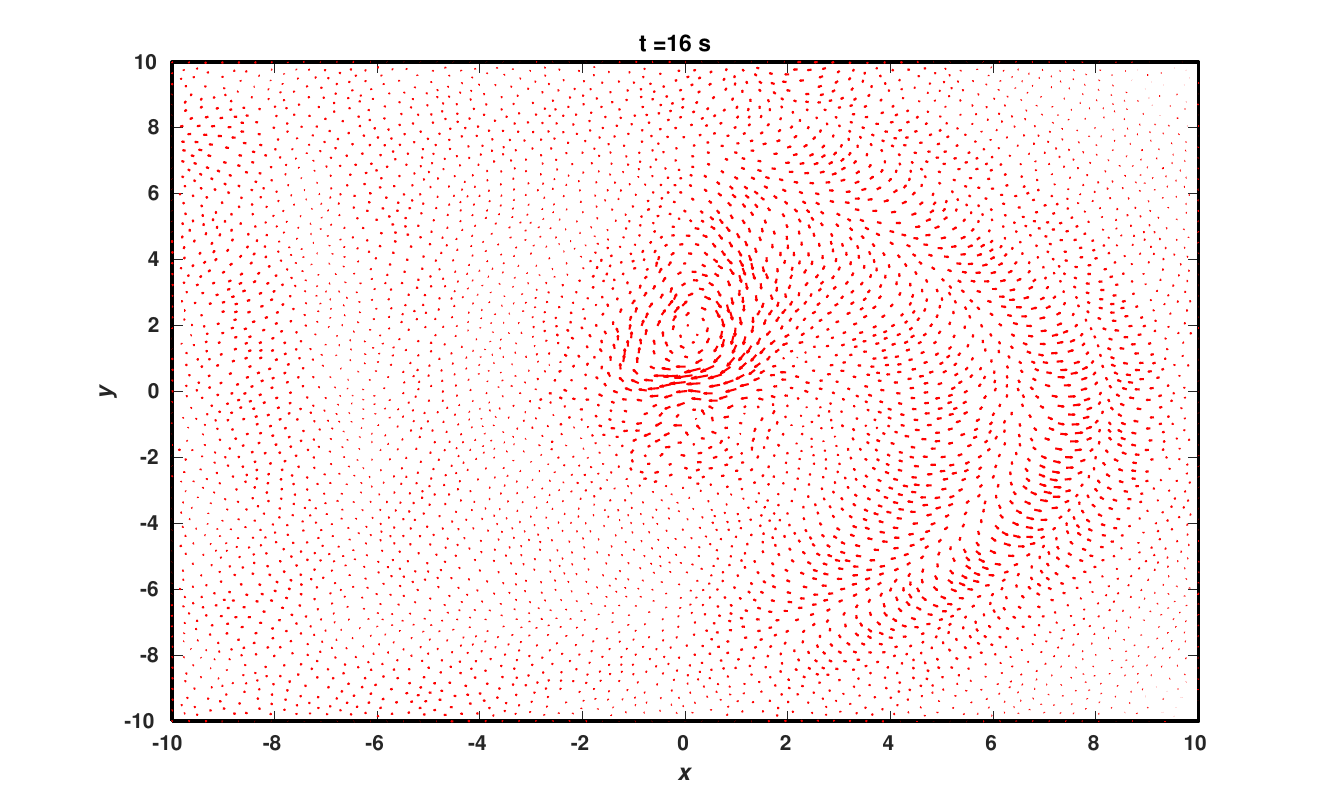} 
	\end{subfigure}
\vspace{.1cm}\\
	\begin{subfigure}[b]{0.33\linewidth}
	\centering
	\includegraphics[width=1.05\textwidth,height=.19\textheight]{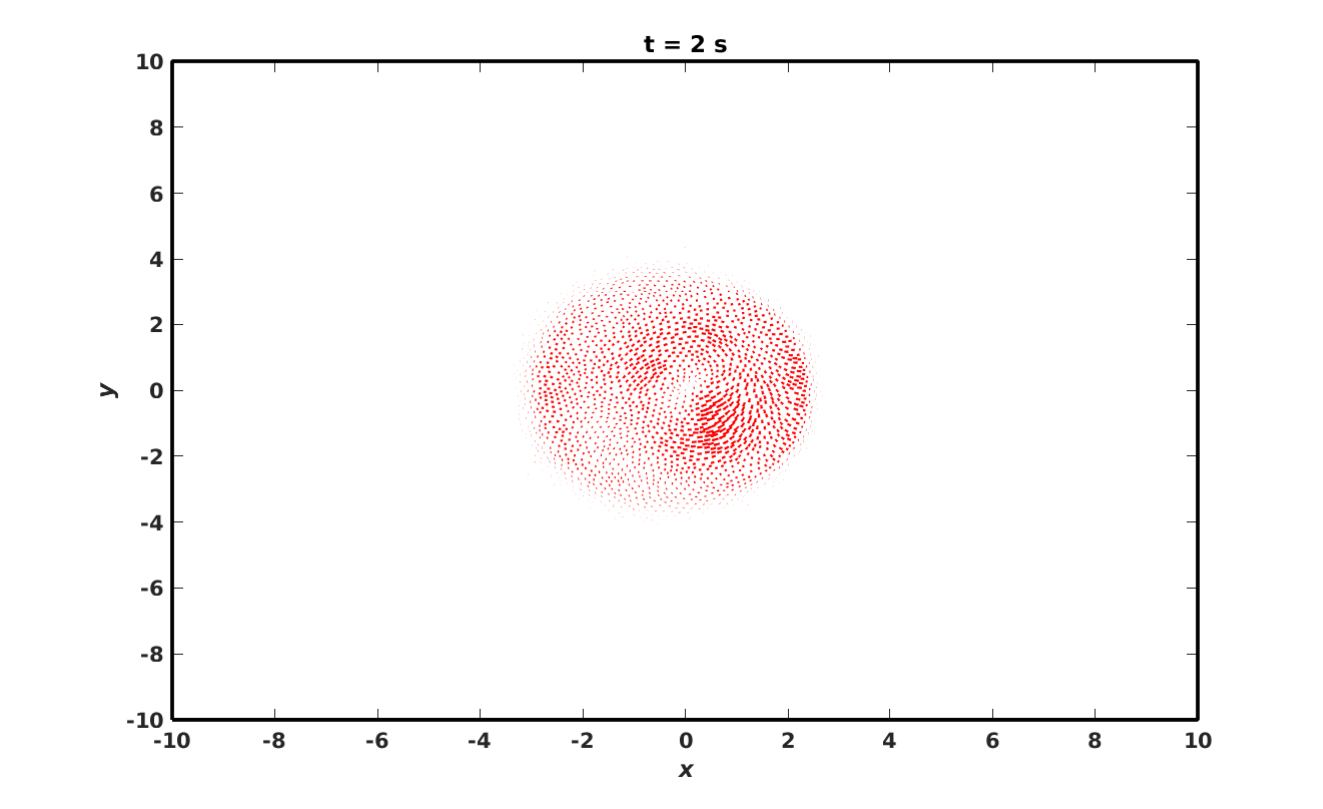} 
\end{subfigure}
\begin{subfigure}[b]{0.33\linewidth}
	\centering
	\includegraphics[width=1.05\textwidth,height=.19\textheight]{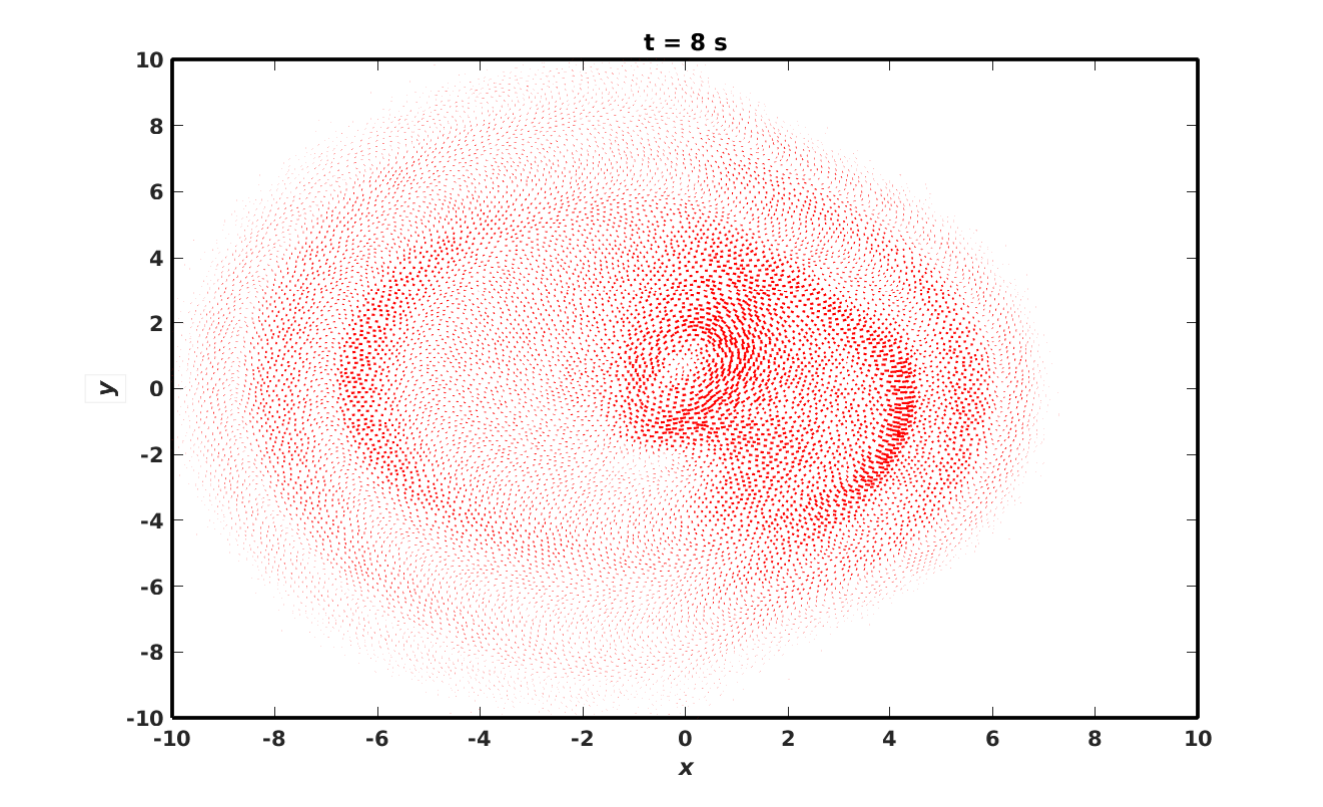} 
\end{subfigure}
\begin{subfigure}[b]{0.33\linewidth}
	\centering
	\includegraphics[width=1.05\textwidth,height=.19\textheight]{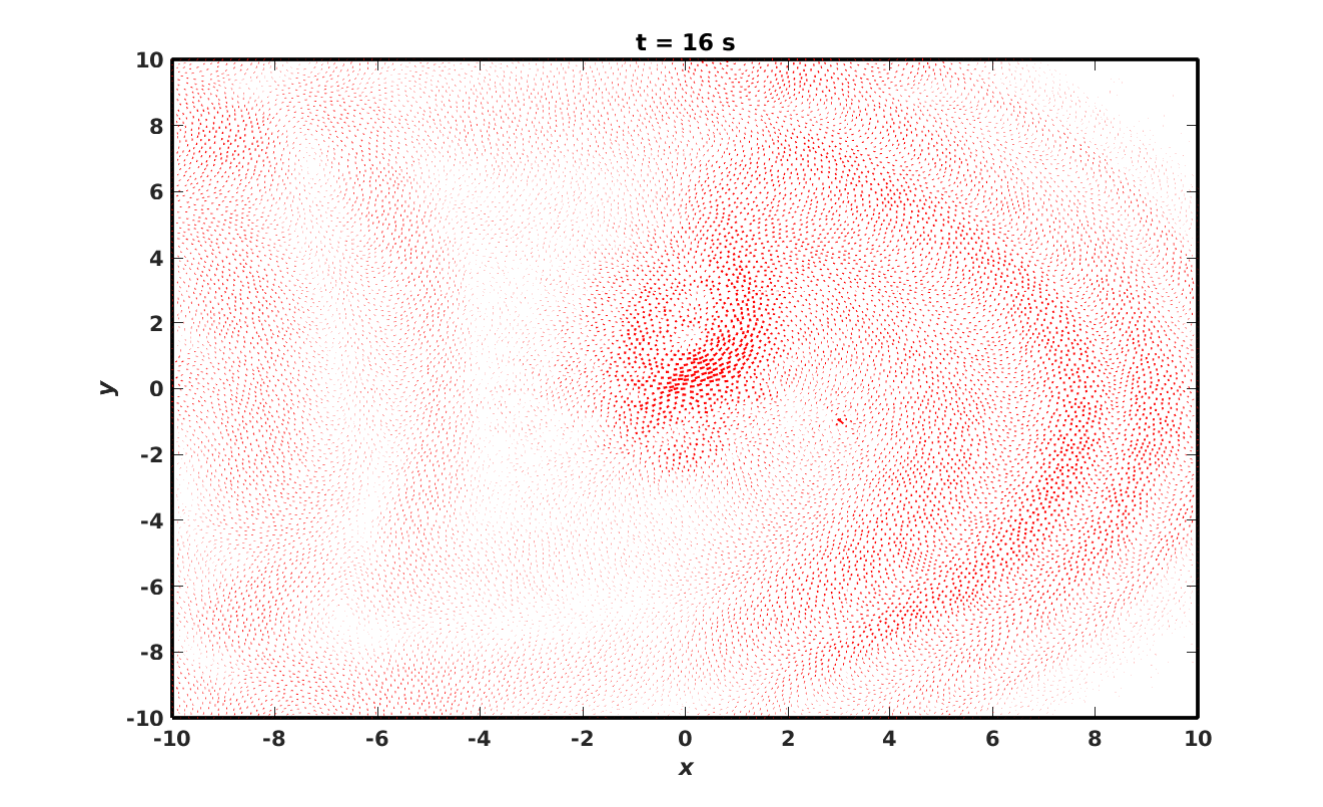} 
\end{subfigure}
	\caption{ Velocity fields corresponding to the plots represented in \figurename{\ref{fig:Circular3}}. First row: mesh with 10040 cells. Second row:  mesh with 40146 cells.
	}
	\label{fig:Circular4}
\end{figure}%

\section{Conclusion and perspectives} \label{sec5}
We have presented in the first part of this paper the extension of the finite volume-characteristics scheme for solving two-dimensional equations of nonlinear conservation laws on unstructured formalism. In the second part, a conservative approach has been presented to approximate the source terms while preserving the stability properties of the homogeneous solver. We note this work was the future goal of many works \cite{benkhaldoun2015family,benkhaldoun2015projection,audusse2014fast,al2019new}.\\
Our approach has several advantages, especially solving steady flows without large numerical errors, thus demonstrating that the proposed scheme preserves the balance related to bathymetry term; it can also compute the numerical flux corresponding to the real state of water flow without relying on Riemann problem solvers.

The method's performance has been evaluated for several test examples; furthermore, in \cite{phongthanapanich2019comparative}, the authors concluded the robustness of the finite volume-characteristics method in other types of equations.
In this work, we have developed a code for calculating free-surface flows over an irregular bottom for complex geometry. This code is based on the solution of the shallow water equations using the FVC scheme in unstructured meshes. This method applies to problems that represent large source terms due mainly to rapid variation and high background irregularity.  Faced with these problems, the methods that use a Riemann approximation solver, which is well suited to the solution of purely hyperbolic equations, often encounter difficulties due to the instability of the numerical solution. In addition to the test cases presented in this paper, there is about fifteen other test cases have been carried out in the context of free-surface flows; each of these tests aims at verifying precisely one or more properties of our code.
It should be noted that some extensions of this approach can be used for modeling realistic applications, notably sediment and pollutant transport. Extending this approximation method to a multi-layer shallow-water model and multi-phase flows \cite{audusse2014fast,brennen2005fundamentals} will also aim for future works. We can also adapt this method for solving the shallow water magnetohydrodynamic equations \cite{qamar2006application,zia2014numerical}.

%
\section*{Acknowledgements} \label{sec6}
The authors thank I. El Mahi, and A. Ratnani for fruitful discussions and helpful comments. This work was partially supported by the HPC Project Alkhwarizmi department, MSDA-UM6P.

\small
\bibliographystyle{abbrv}
\bibliography{./Chapitres/references}
\end{document}